\documentclass{amsart}

\usepackage{amssymb}

\newtheorem{thm}{Theorem}

\newtheorem{lem}[thm]{Lemma}
\newtheorem{cor}[thm]{Corollary}

\newtheorem{prop}[thm]{Proposition}

\newtheorem{conj}[thm]{Conjecture}
   
\theoremstyle{definition}
\newtheorem{defn}[thm]{Definition}

\newtheorem{say}[thm]{}
\newtheorem{exmp}[thm]{Example}

\newtheorem{prob}[thm]{Problem}
\newtheorem{probs}[thm]{Problems}

\newtheorem{rem}[thm]{Remark}          
\newtheorem{note}[thm]{Note}            
\newtheorem{ack}{Acknowledgments}

\newtheorem{defn-thm}[thm]{Definition--Theorem}  

\theoremstyle{remark}

\renewcommand{\c}[0]{{\mathbb C}}  

\renewcommand{\o}[0]{{\mathcal O}} 
\newcommand{\z}[0]{{\mathbb Z}}

\renewcommand{\r}[0]{{\mathbb R}}

\newcommand{\p}[0]{{\mathbb P}}
\newcommand{\f}[0]{{\mathbb F}}
\newcommand{\q}[0]{{\mathbb Q}}

\newcommand{\qtq}[1]{\quad\mbox{#1}\quad}
\newcommand{\spec}[0]{\operatorname{Spec}}
\newcommand{\pic}[0]{\operatorname{Pic}}

\newcommand{\rank}[0]{\operatorname{rank}}
\newcommand{\mult}[0]{\operatorname{mult}}

\newcommand{\supp}[0]{\operatorname{Supp}}    
\newcommand{\red}[0]{\operatorname{red}}    
    
\newcommand{\im}[0]{\operatorname{im}}

\newcommand{\sing}[0]{\operatorname{Sing}}

\newcommand{\tors}[0]{\operatorname{tors}}

\newcommand{\weil}[0]{\operatorname{Weil}}

\newcommand{\ric}[0]{\operatorname{Ricci}}

\newcommand{\rdown}[1]{\lfloor{#1}\rfloor}

\newcommand{\lcm}[0]{\operatorname{lcm}}

\def\into{\DOTSB\lhook\joinrel\to}

\begin{document}
\bibliographystyle{amsalpha}

\title{Einstein metrics on 5-dimensional Seifert bundles}
\author{J\'anos Koll\'ar}
\today
\maketitle

The aim of this paper is to study the existence 
of positive Ricci curvature metrics, and especially
Einstein metrics, on 5-dimensional Seifert bundles.

Seifert fibered 3--manifolds were introduced and studied in \cite{seif}.
These are 3--manifolds $M$ which admit a
differentiable map $f:M\to F$ to a surface $F$
such that every fiber is a circle. 
There are finitely many points  $p_1,\dots,p_n\in F$ 
such that $M$ is a  circle bundle over $F\setminus\{p_1,\dots,p_n\}$,
and the fibers over $p_i$ naturally appear with some multiplicty
$m_i$. 
It is the {\it orbifold  Euler charactersitic}
$$
e^{orb}(F, (p_1,m_1),\dots, (p_n,m_n)):=e(F)-\sum (1-\tfrac1{m_i})
$$
that by and large determines the geometry of $M$. 
For instance, $M$ has spherical geometry  iff
$e^{orb}>0$ and $H_1(M,\q)=0$ (see, for instance,  \cite{scott}).

Higher dimensional
Seifert fibered manifolds were introduced and investigated
in \cite{or-wa}. 
Roughly speaking, these are $(2n+1)$--manifolds $L$ which admit a
differentiable map $f:L\to X$ to a {\em complex} $n$-manifold $X$
such that every fiber is a circle.  Insisting that the base of a
Seifert bundle be a complex manifold seems very artificial from
the topological point of view, but a remarkable result a
Kobayashi \cite{kob} implies that for our purposes this
is necessary. (See (\ref{koba.thm}) for a precise statement.)

The natural setting seems to be to study
Seifert bundles $f:L\to X$, where the base is a complex 
{\it locally cylic orbifold}. That is, locally it looks like
$\c^n/G$ where $G$ is a cyclic group acting linearly.
There is a divisor $\cup D_i\subset X$ such that
$L\to X$ is a circle bundle over $X\setminus \cup D_i$
and there are natural multiplicities $m_i$ are assigned to the
fibers over each $D_i$, see (\ref{seif.descr}). We call 
$\Delta:=\sum (1-\tfrac1{m_i})D_i$ the {\it branch divisor},
and we denote the base orbifold by $(X,\Delta)$.
The key invariant is the {\it orbifold Chern class}
$$
c_1^{orb}(X,\Delta):=c_1(X)-\sum (1-\tfrac1{m_i})[D_i]\in H^2(X,\q).
$$
Another invariant, which enters into the final picture
surprisingly little, is the Chern class of the Seifert bundle
$c_1(L/X)\in H^2(X,\q)$, to be defined in (\ref{chern.class.defn}).

Let $(X,g)$ be a Riemannian manifold and $L\to X$ a circle
bundle over $X$ with a connection. At each point of
$L$ the connection decomposes the tangent space  into a vertical and a 
horizontal piece. By choosing the standard metric on the
circle fibers, we can lift the
metric $g$ to a metric $g_L$ on $L$.
With suitable care, this also works for
Seifert bundles $f:L\to (X,\Delta)$.

\begin{thm}\label{metric.lift.thm} \cite{kob, bg00, bgn03c}
Let $(X,\Delta)$ be a compact, complex
orbifold and $f:L\to (X,\Delta)$ a Seifert bundle. Assume that
$c_1(L/X)$ is a positive rational multiple of 
$c_1^{orb}(X,\Delta)$.

Then, a positive Ricci curvature orbifold K\"ahler metric 
(resp.\ K\"ahler--Einstein metric) on
$(X,\Delta)$ can be lifted  to a positive Ricci curvature metric 
(resp. Einstein metric) on $L$.

Moreover, the lifted metric is also Sasakian.
\end{thm}

It should be noted that
a negative Ricci curvature K\"ahler--Einstein metric on $(X,\Delta)$
does not lift to a negative Ricci curvature Einstein metric on $L$
in any natural way.

There are only a handful of cases when $L\to X$ is a circle bundle and the
above theorem applies, but 
Boyer and Galicki  \cite{bg00, BG01} observed that  many more cases appear 
when $X$ is allowed to be an orbifold.

The method of Boyer and Galicki starts with a complex hypersurface
$0\in Y\subset \c^n$ with an isolated singularity at the origin
which is invariant under a $\c^*$-action
$(z_1,\dots, z_n)\mapsto (\lambda^{w_1}z_1,\dots, \lambda^{w_n}z_n)$.
The intersection of $Y$ with the unit sphere $L:=Y\cap S^{2n-1}(1)$
is called the {\it link} of $0\in Y$. $L$ is a $(2n-3)$-dimensional real
manifold with an $S^1$-action. The quotient 
$X:=L/S^1\cong (Y\setminus \{0\})/\c^*$ is naturally a
complex orbifold.

Methods of complex singularity theory allow one to
identify $L$ as a manifold, and this leads to a large class
of new Einstein metrics on various spaces, including spheres and exotic
spheres \cite{bgk}.

The aim of this paper is to further generalize 
this construction to
arbitrary
 Seifert bundles. The key advantage of this approach is that we can start
with an  orbifold $(X,\Delta)$ and  construct 
 Seifert bundles over $(X,\Delta)$ later. 
This provides substantially greater flexibility, allowing one to explore
the natural scope of the theory.
The construction of
higher dimensional Seifert bundles was considered in \cite{or-wa}
for $X$ smooth; the general case is treated in \cite{k-seif}.
The integral cohomology of 
Seifert bundles is rather subtle in general, but the
5--dimensional case is quite manageable.

The natural questions to be considered 
can be grouped around four  problems:

\begin{probs}{\ }

\begin{enumerate}
\item \label{prob1}
 Describe all manifolds $L$ which have a Seifert bundle sutructure
$f:L\to  (X,\Delta)$  over an orbifold with
positive orbifold Chern class $c_1^{orb}(X,\Delta)$.
\item \label{prob2}
 Given a manifold $L$, describe all Seifert bundle sutructures
$f:L\to  (X,\Delta)$  over an orbifold with
 $c_1^{orb}(X,\Delta)>0$.
\item \label{prob3}
Construct positive Ricci curvature orbifold K\"ahler--Einstein metrics
on orbifolds $(X,\Delta)$.
\item \label{prob4}
Study the place of the
resulting Einstein metrics in the theory of Einstein manifolds.
\end{enumerate}
\end{probs}

It seems rather artificial to
 separate the first two problems, but 
they are quite different in nature. I believe that
in dimension five Problem (\ref{prob1}) is  doable with the present methods,
whereas Problem (\ref{prob2}) seems  hopeless to me.
The reason for this is connected with the complex
geometry of the quotients $X=L/S^1$. These are
{\it Del Pezzo} surfaces with quotient singularities, also
called {\it log Del Pezzo} surfaces.
While  smooth Del Pezzo surfaces have been classified and understood
for more than a century, log Del Pezzo surfaces occur in bewildering abundance
and complexity; see, for instance,  \cite{miyanishi, keel-mck, shokurov}. 

Nonetheless,
as we  see, Seifert bundles over log Del Pezzo surfaces tend to have
 simple topology. Thus many cases are excluded, and for some of the
extreme cases one can get a complete description. Having a
Seifert bundle with simple topology imposes only mild conditions
on a log Del Pezzo surface, and I do not see how to get a
good description. In all likelihood, the hardest is to
describe all Seifert bundle structures on $S^5$.

By Myers' theorem, the fundamental group of a
compact manifold with positive Ricci curvature is finite.
Therefore we concentrate of those cases when $L$ is simply connected.
We see in (\ref{w2=0.lem}) that for  every simply connected
manifold with such a Seifert bundle structure, the  second
Stiefel--Whitney class is zero.
These 5--manifolds are   well understood topologically:

\begin{thm}\label{smale.thm}\cite{smale}
Let $L$ be a  simply connected compact 5--manifold 
with vanishing second
Stiefel--Whitney class. Then:
\begin{enumerate}
\item $L$ 
 is uniquely determined by $H_2(L,\z)$.
\item 
The torsion subgroup of the second homology
is of the form $\tors H_2(L,\z)\cong A+A$ for some finite Abelian group $A$.
\item For any  finite Abelian group $A$ and for any $n\geq 0$
there is a (unique) $L$ with $H_2(L,\z)\cong \z^k+A+A$.
\end{enumerate}
\end{thm}

If $H_2(L,\z)$ is torsionfree of rank $k$ 
then $L$ is the connected sum
of $k$ copies of $S^2\times S^3$.  For any  $k$, 
Einstein metrics on these were constructed
in \cite{bgn02, bgn03, bg03, k-s2s3}. 
 By contrast, the first   result of this paper
shows that
very few of the
possible torsion subgroups do occur.

\begin{thm}\label{h2.tors.thm} Let $L$ be a compact 5--manifold  
such that $H_1(L,\z)=0$.
Assume that $L$ has a Seifert bundle structure
$f:L\to (S,\Delta)$ with $c_1^{orb}(S,\Delta)>0$. 
Then  the torsion subgroup of the second homology, $\tors H_2(L,\z)$,
is one of the following:
\begin{enumerate}
\item $(\z/m)^2$ for any $m$,
\item $(\z/5)^4$ or $(\z/4)^4$,
\item $(\z/3)^4,(\z/3)^6$ or $(\z/3)^8$,
\item $(\z/2)^{2n}$ for any $n$.
\end{enumerate}
Conversely, all these cases do occur, even for manifolds
with Einstein metrics.
\end{thm}

One can be even more precise  if $H_2(L,\z)$ is torsion, that is, 
when $L$ is a {\it rational homology sphere}.
The following characterization uses the notion of the 
{\it orbifold fundamental group} $\pi_1^{orb}(X,\Delta)$
and its abelianization $H_1^{orb}(X,\Delta)$,
to be defined in (\ref{orb.pi1.defn}).

\begin{thm}\label{rhs.seif-strs.thm} 
There is a one--to--one correspondence between
\begin{enumerate}
\item Seifert bundle structures  $f:L\to (S,\Delta)$
on 5-dimensional, compact
rational homology spheres with $H_1(L,\z)=0$, and
\item compact, complex, 2--dimensional,
 locally cyclic orbifolds $(S,\Delta)$ with 
$H_2(S,\q)=\q$ and  $H_1^{orb}(S,\Delta)=0$.
\end{enumerate}
\noindent Under this correspondence,
$\pi_1(L)=\pi_1^{orb}(S,\Delta)$.
\end{thm}

There are many such orbifolds $(S,\Delta)$
 where $c_1^{orb}(S,\Delta)$ is negative,
but very few when $c_1^{orb}(S,\Delta)$ is positive.
Three infinite series  were found in \cite{bg-5dim},
giving  Einstein metrics on rational homology spheres $L$ with 
$H_2(L,\z)=(\z/m)^2$ for any $m$ not divisible by 6.
 We give a complete classification of all cases
when  $H_2(L,\z)$  contains a torsion element of  large enough order.
 Besides the above 3 series, there is only one more
infinite series, but dozens, probably hundreds, of sporadic examples.

\begin{thm}\label{RHS.main.thm}
 Let $L$ be a 5-dimensional compact rational homology sphere
which has a Seifert bundle structure
$f:L\to (S,\Delta)$ with $c_1^{orb}(S,\Delta)>0$.
Assume that  $H_1(L,\z)=0$ and 
 $H_2(L,\z)$  contains a torsion element of  order
 at least 12. Then:
\begin{enumerate}
\item  $H_2(L,\z)=(\z/m)^2$ for some $m$
 not divisible by 30,
\item $L$ is simply connected,
\item the number of Seifert bundle structures
varies between 1 and 4, depending on $m\mod 30$,
\item each Seifert bundle structure gives rise to a 2--parameter family
of Einstein metrics,  naturally
parametrized by the moduli space of  genus 1 curves.
\end{enumerate}
(See (\ref{main.series.h2.tors}) for a  detailed description of
the orbifolds $(S,\Delta)$.)
\end{thm}

\begin{rem}
There are further examples with
 $H_2(L,\z)=(\z/m)^2$ for every $m\leq 11$, see  (\ref{P123.to.11.exmp},
\ref{E8.to.10.exmp}).

 The rank of $H_2(L,\z)$ and its torsion subgroup are not independent
 of each other, but I have only fragmentary results  describing their
interaction. 
\end{rem}

\begin{thm}\label{RHS.main2.thm} Let $L$ be a 5-dimensional compact manifold
which has a Seifert bundle structure
$f:L\to (S,\Delta)$ with $c_1^{orb}(S,\Delta)>0$.
 Assume that $H_1(L,\z)=0$ and
write $H_2(L,\z)=\z^n+A+A$.
\begin{enumerate}
\item  If $H_2(L,\z)=\z^n+(\z/m)^2$ for some $m\geq 12$ then 
\begin{enumerate}
\item $L$ is  simply connected and  $n\leq 8$,
\item there are 93 such cases for every $m\geq 12$,
 (see (\ref{simpconn.cyc.DV.thm}) for the complete list)
\item each Seifert bundle structure yields an at least  2--parameter family of
Einstein metrics.
\end{enumerate}
\item If $H_2(L,\z)=\z^n+(\z/5)^4$ then 
\begin{enumerate}
\item $L$ is a simply connected rational homology sphere,
\item $L$ and the Seifert bundle structure are unique, and 
\item $L$ admits a 4--dimensional family of
Einstein metrics,  naturally
parametrized by the moduli space of  genus 2 curves.
\end{enumerate}
\end{enumerate}
\end{thm}

\begin{say}[K\"ahler--Einstein metrics] In Section 7  we
construct positive Ricci curvature K\"ahler--Einstein 
metrics on certain 2--dimensinal orbifolds.
While these examples are mostly new, the method is
the same as in \cite{jk1}, relying on earlier
works of \cite{nadel, dk}. Thus nothing essentially new
is added to Problem (\ref{prob3}).
\end{say}

\begin{say}[Sasakian  manifolds]
While I prefer to think of Seifert bundles as a
topological object $L$ 
 associated to an algebro--geometric object
$(X,\Delta)$, they have a  natural place
within the framework of Sasakian geometry.
(See  \cite{bg-review} for a recent survey paper.)

Roughly speaking, a {\it quasi--regular Sasakian manifold}
is a  Seifert bundle $L$ 
over a K\"ahler orbifold $(X,\Delta)$  plus a
metric on $L$ which optimally matches the
orbifold K\"ahler metric on $(X,\Delta)$.

In the language of Sasakian geometry, the main results 
of this paper are the
following:
\begin{enumerate}
\item Corollary (\ref{rhs.seif-restriction}) gives
  topological restrictions for a 5--dimensional
rational homology sphere to admit a quasi--regular Sasakian structure.
\item Theorem  (\ref{h2.tors.thm}) shows that 
  most  5--manifolds do not
 admit a positive quasi--regular Sasakian structure.
\item Theorems  (\ref{RHS.main.thm}) and (\ref{RHS.main2.thm}) classify all 
 positive quasi--regular Sasakian structures 
on certain 5--manifolds.
\item We also get new examples of quasi--regular Sasakian--Einstein
metrics, though the examples discovered by Boyer and Galicki
already cover almost all cases allowed by Theorems (\ref{h2.tors.thm}) and
 (\ref{RHS.main.thm}).
\end{enumerate}
\end{say}

\begin{say}[Description of the sections]
Basic results on Seifert bundles and on the Kobayashi construction
are recalled in Section 1. 
In the most general case, we are led to study
Seifert bundles over orbifolds which locally look like
the quotient of $\c^n$ by a cyclic group.
These  are studied in Section 2.

Our main aim is to understand 5--dimensional Seifert bundles
where the base is a log Del Pezzo surface (\ref{lDP.defn}).
These are  rational surfaces with quotient singularities.
In Section 3 we see how to compute the topological (co)homology
of these surfaces in terms of their algebraic geometry. A
key point is to compute everything with $\z$-coefficient.

The cohomology groups of a  Seifert bundle
$f:L\to S$ 
are computed by a  Leray
spectral sequence, and we study it in Section 4.
The spectral sequence degenerates at $E_3$ and
we get a pretty complete description of the $E_2$-term.
The  differential  $E^{0,1}_2\to E^{2,0}_2$ is identified with a 
first Chern class. 
A key observation is that torsion in $H_2(L,\z)$ 
comes from curves $C\subset S$ of genus at least 1  
such that every fiber of $f$ above $C$ is multiple.

Log Del Pezzo surfaces have been the objects of intense investigation
(see \cite{miyanishi, keel-mck, shokurov}). While  smooth Del Pezzo surfaces
form a well understood and easy to describe class,
log Del Pezzo surfaces are rather numerous.
Nonetheless, it is quite rare that  a Seifert bundle  can have 
multiple fibers over a
curve of genus at least 1,
and many of these are classified in Section 5.

Del Pezzo surfaces with cyclic Du Val singularities
and simply connected smooth part are listed in Section 6.

Section 7 establishes the existence of K\"ahler--Einstein metrics on
most of the surfaces considered previously. 

A characterization  of Seifert bundle
structures on rational homology spheres is given in Section 8
in terms of algebraic geometry.  There are probably
too many cases for a meaningful
classification.  
Here we also  collect the details to get
 proofs of the main theorems.

There is a close relationship between manifolds with
Seifert bundles  and
 links of 3--dimensional log terminal
singularities. Some of the resulting open  questions are mentioned
 in Section 9.
\end{say}

\begin{ack} I thank Ch.\ Boyer, K.\ Galicki and  J. McKernan
for many useful conversations and e-mails.
Research
was partially supported by the NSF under grant number
DMS-0200883. 
\end{ack}

\section{Einstein  metrics on Seifert bundles}

Let $(X,g)$ be a Riemannian manifold and $L\to X$ a circle
bundle over $X$ with a connection $D$. At each point of
$L$ the connection decomposes the tangent space  into a vertical and a 
horizontal piece. By choosing the standard metric on the
circle fibers, we can lift the
metric $g$ to a metric $g_D$ on $L$. Kobayashi \cite{kob}
asked:
If $(X,g)$ is an Einstein manifold, when is
$(L,g_D)$ Einstein?

It turns out that this happens only rarely.

\begin{thm}\cite{kob}\label{koba.thm}
 Notation and assumptions as above. Then
$(L,g_D)$ is an Einstein manifold iff and only if one of the
following two conditions are satisfied:
\begin{enumerate}
\item $(X,g)$ is Ricci flat and $D$ is flat.
\item $X$ is a complex manifold, $g$ is the real part
of a K\"ahler--Einstein metric $\omega$, $\ric(\omega)>0$ and the curvature of $D$ is a
positive  multiple of  $\ric(\omega)$.
\end{enumerate}
\end{thm}

The first case is essentially trivial, and here we are
interested in the second case.

In any given dimension, there are only finitely many deformation
types of compact  complex manifolds with positive Ricci curvature
(cf.\ \cite{nadel-fano, campana, kmm3}),
so (\ref{koba.thm}.2) gives only a handful of new  Einstein manifolds.

If $X$ is a compact  complex manifold with positive Ricci curvature
then $H^i(X,\o_X)=0$ for every $i>0$, in particular,
every topological circle bundle over $X$ can be written uniquely
as the unit circle bundle of a holomorphic $\c^*$-bundle
 $f:Y\to X$.
Thus the second alternative of (\ref{koba.thm})
can be formulated completely in terms of complex geometry.
This is the setting where one can generalize the
Kobayashi construction.

\begin{defn} Let $X$ be a normal complex  space.
 A {\it Seifert $\c^*$-bundle}
over $X$ is a normal complex space $Y$
together with a morphism $f:Y\to X$ and a 
$\c^*$-action on $Y$
satisfying the following two conditions.
\begin{enumerate}
\item $f$ is Stein (that is, the preimage of any open  Stein subset of $X$ 
is Stein)
and $\c^*$-equivariant
(with respect to the trivial action on $X$).
\item For every $x\in X$, the $\c^*$-action on the reduced fiber
$Y_x:=\red f^{-1}(x)$, 
$\c^*\times Y_x\to Y_x$ is $\c^*$-equivariantly biholomorphic to
the natural  $\c^*$-action on $\c^*/\mu_{m}$ for some $m=m(x,Y/X)$,
where $\mu_{m}\subset \c^*$ denotes the group of $m$th roots of unity.
\end{enumerate}
\noindent The number $m(x,Y/X)$ is called the {\it multiplicity} of the
Seifert fiber over $x$.
\end{defn}

One can always assume that the $\c^*$-action is effective,
that is, $m(x,Y/X)=1$ for general $x\in X$.

Note that even if $Y$ is smooth, $X$ can have quotient singularities.

A classification of Seifert $\c^*$-bundles for $X$ a smooth
manifold with $H_1(X,\z)$ torsion free is given in
\cite{or-wa}. Many other cases are described in
\cite{dolgachev, pinkham, demazure, fl-za}. The general case
is discussed in detail in \cite{k-seif}. We recall the
relevant facts below.

\begin{say}[Description of Seifert $\c^*$-bundles]
\label{seif.descr}
Let $f:Y\to X$ be a Seifert $\c^*$-bundle. The set of points
$\{ x\in X: m(x,Y/X)>1\}$ is a closed analytic subset of $X$.
It can be written as the union of Weil divisors $\cup D_i$ and
of a  subset of codimension at least 2 contained in $\sing X$.
The latter will not be relevant to us. The multiplicity
$m(x,Y/X)$ is constant on a dense open subset of each $D_i$, this common
value is called the multiplicity of the Seifert $\c^*$-bundle over $D_i$;
denote it by $m_i:=m(D_i)$. We call the $\q$-divisor
 $\Delta:=\sum (1-\tfrac1{m_i})D_i$
the {\it branch divisor} of $f:Y\to X$. We frequently write
$f:Y\to (X,\Delta)$ to indicate the branch divisor.

\end{say}

\begin{thm}\cite[Thm.7]{k-seif}
Let $X$ be a normal complex space with quotient singularities
and $\Delta:=\sum (1-\tfrac1{m_i})D_i$ a $\q$-divisor.
There is a one--to--one correspondence between
Seifert $\c^*$-bundles $f:Y\to (X,\Delta)$ and the following data:
\begin{enumerate}
\item For each $D_i$ an integer $0\leq b_i<m_i$, relatively prime to $m_i$, 
and
\item a linear equivalence class of Weil divisors $B\in \weil(X)$.
(In all the cases that we consider in this paper,
$\weil(X)\cong H_2(X,\z)$.)
\end{enumerate}
\end{thm}

\begin{say}[Algebraic construction] \cite[Thm.7]{k-seif}
 Given the above data, set
$$
\begin{array}{lcl}
S(B, \sum \tfrac{b_i}{m_i}D_i)& :=&\sum_{j\in \z} 
\o_X\bigl(jB+\sum_i \rdown{jb_i/m_i}D_i\bigr)\qtq{and}\\
Y(B, \sum \tfrac{b_i}{m_i}D_i)& :=&\spec_X S(B, \sum \tfrac{b_i}{m_i}D_i),
\end{array}
$$
where $\rdown{x}$ denotes the round down or integral part of a real number.

There is a natural $G_m$-action on $S(L, \sum \tfrac{b_i}{m_i}D_i)$
where the $j$th summand is the  $\lambda^{j}$ eigensubsheaf.
This defines a $\c^*$-action on $Y(L, \sum \tfrac{b_i}{m_i}D_i)$.
\end{say}

\begin{say}[Topological construction]\cite[3.8]{or-wa}
 Set  $M(\Delta)=\lcm(m_i)$. 
If $X$ is smooth then 
$Y\setminus f^{-1}(\supp \Delta)$ can be written as
an $M(\Delta)$-sheeted covering space of the $\c^*$-bundle over
$X\setminus \supp \Delta$ with Chern class
$M(\Delta)\cdot ([B]+ \sum \tfrac{b_i}{m_i}[D_i])$.
The precise description needs the calculation of
$H_1(Y\setminus f^{-1}(\supp \Delta),\z)$ given in \cite[2.1]{or-wa}.
(Note that $M(\Delta)\cdot ([B]+ \sum \tfrac{b_i}{m_i}[D_i])$ can be viewed
as a well defined integral cohomology class, since $m_i| M(\Delta)$.
Dealing with torsion in $H^2(X,\z)$ is, however, delicate, and
\cite{or-wa} assumes that there is no torsion.)

If $X$ is singular, then one first proves that $f:Y\to (X,\Delta)$
is uniquely determined by its restriction to the
smooth locus of $X$ and then the previous method applies.
\end{say}

\begin{defn}\label{chern.class.defn}
 The {\it Chern class} of $f:Y\to (X,\Delta)$
is defined as
$$
c_1(Y/X):=[B]+ \sum \tfrac{b_i}{m_i}[D_i]\in H^2(X,\q).
$$
\end{defn}

\begin{defn}[Orbifolds]\label{orbif.defn}
  An {\it orbifold} is 
 a normal, compact,  complex space $X$
 locally given by charts
written as  quotients of smooth coordinate charts.
That is, $X$ can be covered by
open charts $X=\cup U_i$ and for each $U_i$ there
is a smooth complex space $V_i$ and a finite group
$G_i$ acting on $V_i$ such that $U_i$ is
biholomorphic to the quotient space $V_i/G_i$.
The quotient maps are denoted by $\phi_i:V_i\to U_i$.

The compatibility condition between the charts
that one needs to assume is that  there are
global divisors $D_j\subset X$ and ramification
indices $m_j$ such that $D_{ij}=U_i\cap D_j$
and $m_{ij}=m_j$ (after suitable re-indexing).

It is convenient to codify these data
by a single $\q$-divisor, called the {\it branch divisor}
of the orbifold,
$\Delta:=\sum (1-\tfrac1{m_j})D_j$.

The  {\it orbifold canonical class} is $K_X+\sum (1-\frac1{m_i})D_i$.
This is the negative of the {\it orbifold Chern class}
$$
c_1^{orb}(X,\Delta):=c_1(X)-\sum (1-\tfrac1{m_i})[D_i]\in H^2(X,\q).
$$
\end{defn}

\begin{exmp}\label{seif.gives.orbifold}
 Let $f:Y\to (X,\Delta)$ be a Seifert $\c^*$-bundle
with $Y$ smooth. For $x\in X$ pick any $y\in f^{-1}(x)$
and a $\mu_{m}$-invariant smooth hypersurface $V_x\subset Y$ transversal to
$\red f^{-1}(x)$ for $m=m(x,Y/X)$.
 Then $\{\phi_x:V_x\to U_x:=V_x/\mu_{m}\}$
gives an orbifold structure on $X$. The orbifold branch divisor
coincides with the branch divisor of the Seifert bundle defined in
(\ref{seif.descr}).

Note that the orbifolds coming from a smooth Seifert bundle
have the additional property that each $U_x$ is a quotient
by a cyclic group $\mu_{m}$. Such an orbifold is called
{\it locally cyclic}.

We  usually identify $ (X,\Delta)$ with this orbifold structure.
\end{exmp}

\begin{defn}[Metrics on orbifolds]\label{orb.met.defn}
A {\it Hermitian metric} $h$ on the orbifold $(X,\Delta)$
is  a Hermitian metric $h$ on $X\setminus(\sing X\cup \supp\Delta)$
such that for every chart $\phi_i:V_i\to U_i$
the pull back $\phi_i^*h$ extends to a
Hermitian metric on $V_i$.
One can now talk about curvature, K\"ahler metrics, K\"ahler--Einstein metrics
on orbifolds. 
\end{defn}

\begin{defn} As a real Lie group, $\c^*\cong S^1\times \r$, thus
every Seifert $\c^*$-bundle contains a real hypersurface
$L\subset Y$ with an $S^1$-action.
We call $f:L\to (X,\Delta)$ a {\it Seifert $S^1$-bundle}
or simply a {\it Seifert bundle}. (In dimension 3, these are the original
Seifert bundles.)

$Y$ retracts to $L$, thus they have isomorphic
homology and homotopy groups.
\end{defn}

\begin{thm} \cite{kob, bg00}\label{exst.of.E.metric}
Let $f:Y\to (X,\Delta=\sum (1-\tfrac1{m_i})D_i)$ be a 
Seifert $\c^*$-bundle and $f:L\to (X,\Delta)$ the corresponding
Seifert $S^1$-bundle. Assume that $L$ is a manifold.
Then  $L$ admits an $S^1$-invariant Einstein metric 
with positive Ricci curvature  if and only if the following hold.
\begin{enumerate}
\item The orbifold canonical class
$K_X+\Delta$ is anti ample and 
there is an orbifold  K\"ahler--Einstein metric on $(X,\Delta)$.
\item The Chern class   $c_1(Y/X)$  is a negative multiple of
$K_X+\Delta$. \qed
\end{enumerate}
\end{thm}

From the point of view of algebraic geometry, the
most useful property is the first part of
(\ref{exst.of.E.metric}.1). This class of orbifolds
have their own name.

\begin{defn} An orbifold $(X,\Delta)$ is called {\it Fano}
or {\it log Fano} if the orbifold canonical class
$K_X+\Delta$ is anti ample.
\end{defn}

\section{Smooth Seifert bundles}

\begin{say}[Locally cyclic orbifolds]\label{loc.cyc.orb.say}
Let $f:Y\to (X,\Delta)$ be a Seifert $\c^*$-bundle, $\dim X=n$.
Pick a point $x\in X$ and assume that $Y$ is smooth along
$f^{-1}(x)$. As we saw in (\ref{seif.gives.orbifold}),
$(X,\Delta)$ is then a locally cyclic  orbifold near $x$. 
By diagonalizing the cyclic group action, we see that 
locally $x\in X$ is biholomorphic to  $D^n/\mu_{m}$ for  $m=m(x,Y/X)$
where
$D^n\subset \c^n$ is the $n$-dimensional polydisc
and  the
irreducible components of $\Delta$ passing through $x$
are the quotients of (some of) the coordinate hyperplanes in $D^n$.

A straightforward explicit computation (cf.\ \cite[22--25]{k-seif})
shows that $m=r\cdot m_1\cdots m_n$
where
\begin{enumerate}
\item The $m_1,\dots,m_n$ are the multiplicities of the
irreducible components of $\Delta$ passing through $x$.
(We add the necessary number of 1-s if there are fewer than $n$ 
such components.)
\item The $m_1,\dots,m_n$ are pairwise relatively prime.
\item The quotient $\bar D^n:=D^n/\mu_{m_1\cdots m_n}$ is smooth,
 $D^n/\mu_{m}\cong \bar D^n/\mu_{r}$
and the $\mu_{r}$-action is fixed point free outside a codimension 2 subset.
\item Thus $\z/r$ is the local fundamental group of
$X\setminus\sing X$ at $x$, hence
$m(x,Y/X)$  depends only on $(X,\Delta)$ and not on $Y$.
\end{enumerate}
\end{say}

\begin{defn} Let $(X,\Delta)$ be an orbifold given by the charts 
$ \phi_i:V_i\to U_i$.
A divisor $D\subset X$ 
is called {\it orbismooth}
if the set theoretic preimages $\phi_i^{-1}(D)\subset V_i$ are all smooth.

Thus $D$ itself is an orbifold with the induced orbifold structure.

If $\dim X=2$ then $D$ is a curve and
 orbismooth implies smooth. However,
 not all smooth curves in $X$ are orbismooth.

For instance, act on $\c^2$ by $\mu_5$ as 
$(x,y)\mapsto (\epsilon^2x,\epsilon^3y)$.
Then $(x^3-y^2=0)/\mu_5\subset \c^2/\mu_5$ is smooth but
not orbismooth.
\end{defn}

\begin{defn} \label{local.trosion.measures}
Let $(X,\Delta=\sum (1-\frac1{m_i})D_i)$
 be a locally cyclic orbifold.
Such an orbifold behaves very much like a manifold if we use
$\q$-coefficients, but torsion questions become quite
delicate when working with integral (co)homology.
We need several  ways to measure
contributions of the orbifold points.

As above,  for every $x\in X$ we can write the orbifold structure of 
$(X,\Delta)$ in a  suitable neighborhood  as
$D^n/\mu_{m(x)}$ and  the orbifold structure of 
$(X,\emptyset)$ as $\bar D^n/\mu_{r(x)}$.

\begin{enumerate}
\item Set 
$M(x,\Delta):=\lcm(m_i:x\in D_i)$ and $M(\Delta):=\lcm(m(x,\Delta):x\in X)
=\lcm(m_i)$.
\item Set 
$M(x,X)=r(x)$ and $M(X):=\lcm(M(x,X)):x\in X)$
\item Set 
$M(x,X,\Delta)=m(x)$ and $M(X,\Delta):=\lcm(M(x,X,\Delta):x\in X)
=\lcm(m(x):x\in X)$.
\end{enumerate}
As noted in (\ref{loc.cyc.orb.say}),
$M(x,X,\Delta)=M(x,\Delta)\cdot M(x,X)$
and $M(X,\Delta)\vert M(\Delta)\cdot M(X)$ but the latter
can be different.
\end{defn}

\begin{defn} Let $\mu_r$ act on $D^n$
such that the action is fixed point free outside a codimension 2 subset.
Then $\weil(D^n/\mu_r)\cong \z/r$, noncanonically (cf. \cite[24]{k-seif}).
This is called the {\it local class group} of $X$ at $x$, denoted by
$\weil(x,X)$.
One can also
identify this group as the second cohomology of the
smooth part of the quotient $D^n/\mu_r$. Thus if $X$ is an
orbifold, for every singular point of $X$ we
obtain a map $R_x:\weil(X)\to \weil(x,X)$ which is nonzero on
 a Weil divisor $A$ iff $A$ is not Cartier at $x$. Thus
$$
\pic(X)=\cap_{x\in X} \ker R_x\subset \weil(X).
$$
Topologically, we can see this as
$$
H^2(X,\z)=\cap_x \ker [H_{2n-2}(X,\z)\to \weil(x,X)],
$$
where we identify $H^2(X,\z)$ with its image in
$H_{2n-2}(X,\z)$ by capping with the fundamental class.

In particular, 
we can view multiplication by $M(X)$ as a map
$$
M(X):\weil(X)\to \pic(X)\qtq{or} M(X): H_{2n-2}(X,\z)\to H^2(X,\z).
$$
\end{defn}

\begin{defn} 
Let $f:Y\to (X,\Delta)$ be a Seifert $\c^*$-bundle
over the orbifold $(X,\Delta)$. Its
 Chern class 
$$
c_1(Y/X)=[B]+\sum \tfrac{b_i}{m_i}[D_i].
$$
is an element of $H^2(X,\q)$, but certain multiples of it
can be viewed as  well defined integral classes, even in
the presence of torsion. We need three versions of this:
\begin{enumerate}
\item 
$M(\Delta)\cdot c_1(Y/X):=
M(\Delta)\cdot [B]+\sum \tfrac{M(\Delta)}{m_i}b_i\cdot [D_i]$
is a  well defined  element of $\weil(X)$ or of $H_{2n-2}(X,\z)$.
\item $M(x,\Delta)\cdot c_1(Y/X):=
M(x,\Delta)M\cdot [B]+\sum_{i:x\in D_i} \tfrac{M(x,\Delta)}{m_i}b_i\cdot [D_i]$
is a  well defined  element of $\weil(x,X)$.
\item  $M(X,\Delta)\cdot c_1(Y/X):=
M(X,\Delta)\cdot [B]+\sum \tfrac{M(X,\Delta)}{m_i}b_i\cdot [D_i]$
is a  well defined  element of $\pic(X)$ or of $H^2(X,\z)$. 
\end{enumerate}
\end{defn}

In general, a Seifert $\c^*$-bundle over an orbifold is
singular, but the smooth ones are easy to determine:

\begin{prop}[Smoothness criterion]\cite[29]{k-seif}\label{smoothness.crit}
Let $Y:=Y(B,\sum\frac{b_i}{m_i}D_i)$ be a Seifert bundle over
the orbifold $(X,\sum (1-\frac1{m_i})D_i)$. Then $Y$ is smooth
along $f^{-1}(x)$ iff $M(x,\Delta)\cdot c_1(Y/X)$
is a generator of the local class group $\weil(x,X)$.
\qed
\end{prop}

\section{The topology of singular surfaces}

Let $L\to (S,\Delta)$ be a 
Seifert bundle over a compact complex surface.
 As we saw in (\ref{exst.of.E.metric}), the Kobayashi construction
yields an Einstein metric on L only if $-(K_S+\Delta)$ is ample.
In particular, $S$ is always a projective  algebraic surface.
If $L$ is smooth then  $S$ has only quotient singularities
by (\ref{seif.gives.orbifold}). Such  surface
 is always rational (this is an easy special case of  \cite{sakai}
or of \cite{mi-mo}).

In practice, one can understand the
algebraic curves and their intersection theory
on any (singular) rational surface.
The aim of this section is to describe their
topology, especially various (co)homology groups,  in terms of 
algebraic curves.

Let $S$ be a normal compact surface. Let $S^0$
denote the smooth locus and $P_i\in S$ the set of singular
points.
Topologically, near any singular point $S$ is a cone $C_i$ over
a 3--manifold $M_i$ called the {\it link}. 
While we are mainly interested in surfaces with cyclic quotient
singularities, it is equally easy to work with arbitrary
rational singularities.

\begin{defn} Let $0\in F$ be a normal surface singularity
and $g:F'\to F$ a resolution. The singularity $0\in F$ is
{\it rational} if $R^1g_*\o_{F'}=0$. This is independent of the
resolution chosen.

If $0\in F$ is rational then $g^{-1}(0)\subset F'$ is a tree
of smooth rational curves. This implies that  $R^1g_*\z=0$
and $H^1(M,\z)=0$ where $M$ is the link of $0\in F$.
See \cite{mumf-surf} for these and more information on the topology of 
surface singularities.
\end{defn}

\begin{prop}\label{rtl.surf.coh.prop} 
Let $S$ be a normal, compact surface with
rational singularities $P_i$ and links $M_i$.
Assume that $H^1(S,\z)=0$. 
Let $S^0\subset S$ be the smooth locus. Then
\begin{enumerate}
\item 
$H^3(S,\z)=H_1(S^0,\z)$ is torsion. 
\item  $H_2(S^0,\z)\cong H^2(S,\z)$ is torsion free.
\item $H^2(S^0,\z)\cong H_2(S,\z)$.
\item We can write $H_2(S,\z)$ as the direct sum 
of a free group $H_2(S,\z)_{tf}$  and a torsion group
isomorphic to $H_1(S^0,\z)$.
\item There is an exact sequence
$$
0\to H^2(S,\z)\to H_2(S,\z)_{tf}+ H_1(S^0,\z)
\to \sum_i H^2(M_i,\z)\to H_1(S^0,\z)\to 0.
$$
\item $H^3(S^0,\z)\cong \z^s$
where  $s$ is the number of singular points.
\end{enumerate}
\end{prop}

Proof. 
Set $S^*:=S\setminus\cup_j\operatorname{interior}(C_j)$.
The inclusion  $S^*\subset S^0$ is a  homotopy equivalence.

The long exact cohomology sequence of the pair
$(S,S^*)$ gives
$$
H^{i}(S,S^*,\z)\to H^i(S,\z)\to H^i(S^*,\z)\to H^{i+1}(S,S^*,\z)\to \cdots
$$
By excision, $H^{i+1}(S,S^*,\z)$ 
is the direct sum of the local terms $H^{i+1}(C_j,M_j,\z)$.
$C_j$ is contractible, so  $H^i(C_j,\z)=0$ for $i\geq 1$,
thus $H^{i+1}(C_j,M_j,\z)=H^i(M_j,\z)$ for $i\geq 1$.

For any rational 
surface singularity  $H^1(M_i,\z)=0$ and $H^2(M_i,\z)$ is torsion,
thus we get an exact sequence
$$
0\to H^2(S,\z)\to H^2(S^*,\z)\to \sum_j H^2(M_i,\z)
 \to   H^3(S,\z)\to H^3(S^*,\z).
$$

Assume next that $H_1(S,\z)=0$.
Let $g:S'\to S$ be a resolution of singularities.
Since $S$ has rational singularities, the fibers of $g$ are all simply 
connected, thus $H_1(S',\z)=0$ and 
hence  $H^3(S',\z)=0$.
Therefore $H^3(S,\z)\to H^3(S^*,\z)$ is the zero map since it 
factors as $H^3(S,\z)\to H^3(S',\z)\to H^3(S^*,\z)$. 
In particular, $H^3(S^*,\z)\cong \sum_j H^3(M_i,\z)\cong \z^s$
if we have $s$ singular points.

Alexander duality for  $S'\supset \cup_j f^{-1}(C_j)$
gives  isomorphisms
$$
H^i(S^*,\z) \cong H_i( S',\cup_j f^{-1}(C_j),\z)=
H_i( S',\cup_j f^{-1}(P_j),\z).
$$

By excision, 
$H_i( S',\cup_j f^{-1}(C_j),\z)=H_i( S,\cup_j C_j,\z)$ which equals 
$H_i(S,\z)$ for $i\geq 2$.
Finally, by the universal coefficient theorem and by
Alexander duality,
\begin{eqnarray*}
\tors H_1(S^*,\z)&\cong &\tors H^2(S^*,\z)\cong \tors 
H_2( S',\cup_j f^{-1}(C_j),\z)\\
&=&\tors H_2(S,\z)=H^3(S,\z).
\end{eqnarray*}
This completes the proof.\qed

\begin{note} The  appearance of $H_1(S^0,\z)$ in (5) at two places
puts severe restrictions on the group $H_1(S^0,\z)$, especially
since $H^2(S,\z)$ is torsion free. Indeed, this implies that
the sequence
$$
0\to  \left(H_2(S,\z)_{tf}/H^2(S,\z)\right)+ H_1(S^0,\z)
\to \sum_i H^2(M_i,\z)\to H_1(S^0,\z)\to 0
$$
is exact. Thus if $\prod_i |H^2(M_i,\z)|$ is square free, we can conclude
right away that $H_1(S^0,\z)=0$.
\end{note}

\begin{cor}\label{ratsurf.simply.conn.cor} 
  Let $S$ be a normal, compact surface with
rational singularities $P_i$ with links $M_i$.
Assume that $H_1(S,\z)=0$.
The following are equivalent
\begin{enumerate}
\item $H_1(S^0,\z)=0$.
\item $|H_2(S,\z)/H^2(S,\z)|=\prod_i |H^2(M_i,\z)|$.
\item The determinant of the intersection matrix on $H^2(S,\z)$
is $\prod_i |H^2(M_i,\z)|$.
\end{enumerate}
\end{cor}

Proof. $H_1(S^0,\z)=0$ iff $H^3(S,\z)=0$ which is equivalent
to (2) by the sequence (\ref{rtl.surf.coh.prop}.6).
Since the pairing 
$H_2(S,\z)_{tf}\times H^2(S,\z)\to \z$ is perfect,
$|H_2(S,\z)/H^2(S,\z)|$ is the same as the 
determinant of the intersection matrix on $H^2(S,\z)$.\qed
\medskip

\begin{defn}\label{orb.pi1.defn}\cite{thur} 
For an orbifold $(X,\Delta)$ the {\it orbifold fundamental group} 
$\pi^{orb}_1(X,\Delta=\sum_i (1-\frac1{m_i})D_i)$
is the fundamental group of $X\setminus(\sing X\cup \supp\Delta)$
modulo the relations: if $\gamma$ is any 
 small loop around $D_i$ at a smooth point
 then $\gamma^{m_i}=1$.

The abelianization of $\pi^{orb}_1(X,\Delta)$
is  denoted by  $H^{orb}_1(X,\Delta)$, called the
{\it abelian orbifold fundamental group.}
\end{defn}

\begin{prop}\label{ratsurf.orbif.simpconn.prop}
 Let $S$ be a normal, projective  surface with
rational  singularities $P_i$ with links $M_i$.
Then $H_1^{orb}(S,\sum_j (1-\frac1{m_j})D_j)=0$ iff
\begin{enumerate}
\item $H_1(S^0,\z)=0$, and
\item   The map $H^2(S,\z)\to \sum_j \z/m_j$ given by
$L\mapsto (L\cdot D_j)\mod m_j$
 is surjective.
\end{enumerate}
\end{prop}

Proof.  By \cite[4.6]{or-wa}, if $H_1(S^0,\z)=0$ then
$H_1^{orb}(S,\Delta)$ is 
 given by generators $g_1,\dots,g_n$ and relations
\begin{enumerate}\setcounter{enumi}{2}
\item $m_jg_j=0$ for $j=1,\dots,n$, and
\item $\sum g_j([D_j]\cap \eta)=0$
for every $\eta\in H_2(X,\z)$.
\end{enumerate}
Thus there is an exact sequence
$$
H_2(S^0,\z)\stackrel{\sigma}{\to} \sum_j \z/m_j\to H_1^{orb}(S,\Delta)\to 0
$$
where $\sigma(\eta)=\sum ([D_j]\cap \eta)g_j$.

By (\ref{rtl.surf.coh.prop}.2),
$H_2(S^0,\z)$ is isomorphic to $H^2(S,\z)$ and the map
$\sigma$
is identified with  taking the intersection number with 
each $D_i$ (modulo $m_i$).
\qed

\begin{rem}\label{p-covers.cor} A very easy to apply special case of
(\ref{ratsurf.orbif.simpconn.prop}) is the following:

If we can write 
$$
\textstyle{\sum_{i: p\vert m_i}}
 a_i D_i\sim  p\cdot(\mbox{integral Weil divisor})
$$
and not all the $a_i$ are divisible by $p$,
then $H_1^{orb}(S,\Delta)\neq 0$.
Indeed, for any $L\in H^2(S,\z)$ this would
give
$$
\textstyle{\sum_{i: p\vert m_i}}
 a_i (L\cdot D_i)\equiv 0\mod p.
$$
\end{rem}

We now start to connect the previous results with the
algebraic geometry of the surface $S$.

\begin{prop} \label{pic.weil.prop}
Let $S$ be a normal compact surface with rational singularities
such that $H^1(S,\o_S)=0$.
 Let $g:S'\to S$
be a resolution of singularities. 
Then
\begin{enumerate}
\item There are injections with torsion free cokernels
$$
\weil(S)\into H_{2}(S,\z)\qtq{and}
\pic(S)\into H^2(S,\z). 
$$
\item $H_{2}(S,\z)/\weil(S)\cong H_{2}(S',\z)/\weil(S')$.
\item If $H^2(S,\o_S)=0$ then both of the injections are isomorphisms.
\end{enumerate}
\end{prop}

Proof. 
The long exact cohomology sequence of the exponential sequence
$$
0\to \z\stackrel{2\pi i}{\to} \o_S\stackrel{exp}{\to} \o^*_S\to 1
$$
gives $\pic(S)\into H^2(S,\z)$ if  $H^1(S,\o_S)=0$ and
this is an isomorphism if $H^2(S,\o_S)=0$ also holds.
The cokernel is torsion free since  $H^2(S,\o_S)$ is torsion free.

Set $A:=\sing S$ and $A':=f^{-1}(A)$.
The homology sequences of the pairs $(S',A')$ and $(S,A)$ give a 
commutative diagram 
$$
\begin{array}{ccccccc}
H_{2}(A',\z)& \to & H_{2}(S',\z)& \to & H_{2}(S',A',\z)& \to & 
H_{1}(A',\z)\\ 
\downarrow 0 && \downarrow f_* && \downarrow\cong && \downarrow 0\\
H_{2}(A,\z)& \to & H_{2}(S,\z)& \to & H_{2}(S,A,\z)& \to & 
H_{1}(A,\z)\\ 
\end{array}
$$
The horizontal maps at the end are 0 since $\dim A=0$.
 $H_{1}(A',\z)=0$ since the resolution graph of a rational
surface singularity is a tree of rational curves. 
This implies that 
$H_{2}(S,\z)\cong H_{2}(S',\z)/H_{2}(A',\z)$.
$H_{2}(S',\z)/\weil(S')\cong H^{2}(S',\z)/\pic(S')$ is torsion free and 
$\im [H_{2}(A',\z)\to H_{2}(S',\z)]\subset \weil(S')$.
Thus
$$
H_{2}(S,\z)/\weil(S)\cong H_{2}(S',\z)/\weil(S')
$$
is torsion free and it is zero iff $H^2(S,\o_S)=0$.\qed

\begin{rem} In general the
natural map  $H^{2}(S,\z)/\pic(S)\to H^{2}(S',\z)/\pic(S')$
is not an isomorphism. As an example, let $S'\subset \p^3$
be a general quartic surface containing a plane conic $C$.
$C\subset S'$ is a $-2$-curve and we contract it to get $S$.
Then $ H^{2}(S,\z)=\ker[H^{2}(S',\z)\to H^2(C,\z)\cong \z]$
and $\pic(S)=\ker[\pic (S')\to \pic(C)\cong \z]$. The first of these
is surjective while the image of the second is $2\z$.
\end{rem} 

\begin{prop} \label{q-planes.prop}
Let $S$ be a normal compact surface with rational singularities.
\begin{enumerate}
\item If  $H^1(S,\q)=0$  then $H^1(S,\o_S)=0$.
\item If, in addition,   $H^2(S,\q)\cong \q$  then $H^2(S,\o_S)=0$.
\end{enumerate}
\end{prop}

Proof. Let $f:S'\to S$ be a resolution, 
then $R^1f_*\z=0$ and $R^1f_*\o_{S'}=0$ since
$S$ has rational singularities. If  $H^1(S,\q)=0$, the first implies that
$H^1(S',\q)=0$ and so  $H^1(S',\o_S)=0$ by Hodge theory
 (cf.\ \cite[Sec.0.7]{gh}).
This gives $H^1(S,\o_S)=0$.

If $H_2(S,\q)\cong \q$ then it is generated by
algebraic cycles. By (\ref{pic.weil.prop}) this implies that
 $H_2(S',\q)\cong \q$  is generated by
algebraic cycles, hence $H^2(S',\o_{S'})=0$ which gives that
$H^2(S,\o_S)=0$.\qed

\begin{say}[Surfaces with  $H^1(F,\o_F)=H^2(F,\o_F)=0$]

By  classification theory, smooth projective surfaces
 with  $H^1(F,\o_F)=H^2(F,\o_F)=0$ fall in three groups.
(See, for instance, \cite[Chap.VI]{bhpv} for a comprehensive treament.)
\begin{enumerate}
\item $F$ is a rational surface,
\item $F$ is an Enriques surface \cite[VIII.15--21]{bhpv}, or
\item $F$ is of general type with $q=p_g=0$ \cite[VII.10]{bhpv}.
\end{enumerate}

We are specially interested in cases when $F$ is obtained
as a resolution of a surface $S$ with rational singularities
and $H^2(S,\q)\cong \q$. There are many such examples
where $F$ is rational but probably very few examples 
where $F$ is Enriques  or of general type.

A nice class of examples is given by the so called 
{\it fake projective planes}, smooth surfaces $F$ with
 $H^1(F,\z)=0$ and $H^2(F,\z)\cong \z$
(cf. \cite[V.1]{bhpv}). For these we can take $S=F$. 
All of these are quotients of the complex unit ball
(this is a rather difficult result of \cite{yau}) hence
have a large fundamental group.

Most surfaces with $q=p_g=0$  are not simply connected
(cf. \cite[VII.10]{bhpv}). The only known simply connected examples have
large Picard number.

With no basis whatsoever, other than the lack of such examples,
I suggest that this may be a general phenomenon.
The following precise form is the one needed in
the classification of Seifert structures on
simply connected rational homology spheres
to be studied in Section 7.

\begin{conj}\label{simpconn.implies.rtl.conj}
 Let $S$ be a projective surface with
quotient singularities such that $H^2(S,\q)\cong\q$ and
$S^0$ is simply connected. Then $S$ is rational.
\end{conj}

\end{say}

For surfaces with $H_2(S,\q)=\q$, one can 
make (\ref{ratsurf.simply.conn.cor})  even more explicit.

\begin{cor}\label{ratsurf.rho=1simply.conn.cor} 
Let $S$ be a normal, projective surface with
rational  singularities $P_i$ with links $M_i$. Assume that
$H_1(S,\z)=0$ and $H_2(S,\q)=\q$.
 Then
the following  conditions are equivalent
\begin{enumerate}
\item $H_1(S^0,\z)=0$.
\item $\weil(S)\cong \z$.
\item Each $H^2(M_i,\z)$ is cyclic, their orders $m_i$ are pairwise 
relatively prime and there is a Weil divisor $E$ which generates
the  $H^2(M_i,\z)$ for every $i$.
\item There is a Weil divisor $D$ with $(D^2)=1/\prod m_i$.
\item There is a Cartier divsior $H$ and a Weil divisor $D$ such that
$(D\cdot H)=1$ and $(H^2)=\prod m_i$.
\end{enumerate}
\end{cor}

Proof. 
By (\ref{q-planes.prop}) and (\ref{pic.weil.prop}), $\weil(S)\cong H_2(S,\z)$.
Thus
the equivalence of the first 3 are clear from the
sequence (\ref{rtl.surf.coh.prop}.6). The last two conditions are
reformulations of (\ref{ratsurf.simply.conn.cor}.2).\qed 
\medskip

The vanishing of the abelian orbifold fundamental group
is also given by a simple algebraic condition:

\begin{cor}\label{rho=1.orbH1=0.cor}
 Let $S$ be a normal, projective, rational surface with
rational  singularities $P_i$ with links $M_i$.  Assume that
$H_1(S,\z)=0$ and $H_2(S,\q)=\q$. 
Then $H_1^{orb}(S,\sum_j (1-\frac1{m_j})D_j)=0$ iff
\begin{enumerate}
\item $\weil(S)\cong \z$, 
\item the $m_j$ are pairwise relatively prime, and
\item  $m_j$ is relatively prime to the degree
$\deg D_j\in\z\cong \weil(S)$  for every $j$.\qed
\end{enumerate}
\end{cor}

\section{The topology of 5-dimensional Seifert bundles}

The aim of this section is to describe the
(topological) cohomology groups of a  Seifert bundle
over a rational surface in terms of  invariants
of the base orbifold.

The following is a straightforward generalization of the
computation of the fundamental group of 3--dimensional Seifert 
bundles, cf.\ \cite{seif}.

\begin{prop}\label{fund.gr.sequence}\cite[5.7]{HaSa91}, \cite[Prop.50]{k-seif}
 Let $f:L\to (X,\Delta)$ be a Seifert bundle
such  that $L$ is smooth. There is
an exact sequence
$$
\pi_1(\c^*)\to \pi_1(L)\to \pi^{orb}_1(X,\Delta)\to 1.\qed
$$
\end{prop}

The determination of $\pi_1(L)$
seems rather tricky in general, but  its abelianization
 is
fully computable.
The case when $H_1(X^0,\z)=0$ is especially easy to state.

\begin{prop}\cite[4.6]{or-wa}\label{OW.h1.lem}
 Let $X$ be a complex manifold such that $H_1(X,\z)=0$ and let 
$D_1,\dots,D_n\subset X$ be smooth divisors intersecting transversally.
 For any divisor $B$, $H_1(Y(B,\sum\frac{b_i}{m_i}D_i),\z)$
is given by generators $k,g_1,\dots,g_n$ and relations
\begin{enumerate}
\item $m_ig_i+b_ik=0$ for $i=1,\dots,n$, and
\item $k([B]\cap \eta)-\sum g_i([D_i]\cap \eta)=0$
for every $\eta\in H_2(X,\z)$.\qed
\end{enumerate}
\end{prop}

\begin{rem} The assumptions imposed on $X$ in
(\ref{OW.h1.lem}) seem somewhat restrictive, but
the result can be used to compute the
first homology of any smooth Seifert $\c^*$-bundle.
Indeed, if $f:Y\to (X,\Delta)$ is any Seifert $\c^*$-bundle,
then let $Z\subset X$ denote the union of all singular points of $X$
and of all singular points of $\cup D_i$. $Z$ has codimension at least 2
in $X$, so $\pi_1(Y)=\pi_1(Y\setminus f^{-1}(Z))$. Thus
$H_1(Y,\z)=H_1(Y\setminus f^{-1}(Z),\z)$
can be computed by applying
(\ref{OW.h1.lem}) to $X\setminus Z$.
\end{rem}

\begin{cor} \label{fiber.trison.coeff}
With notation and assumptions as in (\ref{OW.h1.lem}), set
 $Y:=Y(B,\sum\frac{b_i}{m_i}D_i)$. 
As noted in (\ref{local.trosion.measures}),
 $M(\Delta)\cdot c_1(Y/X)\in H^2(X,\z)$
is well defined and thus can be written as
 $M(\Delta)\cdot c_1(Y/X)=d(Y)\cdot U$
where $d(Y)\in \z$ and $U\in H^2(X,\z)$ is primitive.

Then $d(Y)\cdot k=0$ in $H_1(Y,\z)$.
\end{cor}

Proof. Mutiply (\ref{OW.h1.lem}.2) by $M(\Delta)$ and use that
$$
M(\Delta)g_i=\tfrac{M(\Delta)}{m_i}(m_ig_i)=-\tfrac{M(\Delta)}{m_i}b_ik
$$
to rewrite it as
$$
\bigl((M(\Delta)\cdot c_1(Y/X))\cap \eta\bigr)\cdot k=0
\qtq{for every $\eta\in H^2(X,\z)$.}
$$
Since $M(\Delta)\cdot c_1(Y/X)=d(Y)\cdot (\mbox{primitive class})$,
there is an $\eta$ such that $(M(\Delta)\cdot c_1(Y/X))\cap \eta=d(Y)$.
Thus $d(Y)\cdot k=0$.\qed
\medskip

The following is a very convenient property
of orbifolds with $H^{orb}_1(X,\Delta)=0$.

\begin{prop}\cite[53]{k-seif} \label{c_1.gives.Sb.for.simpconn}
Yet $(X,\Delta)$ be an orbifold and assume that
$H^{orb}_1(X,\Delta)=0$. Then a Seifert $\c^*$-bundle
$f:Y\to (X,\Delta)$ is uniquely determined by its Chern class
$c_1(Y/X)\in H^2(X,\q)$. \qed
\end{prop}

The main result of this section is the following:

\begin{thm} \label{simpconn.over.ratsurf.thm}
Let $f:L^5\to (S,\Delta=\sum_i (1-\tfrac1{m_i})D_i)$
 be a smooth Seifert bundle over a projective surface 
with rational singularities. Assume that $H_1(L,\q)=0$ and   
$H^{orb}_1(S,\Delta)=0$. 
Set $s=\rank H^2(S,\q)$.
\begin{enumerate}
\item The cohomology groups $H^i(L,\z)$ are
$$
\begin{array}{cccccc}
H^0 &H^1 &H^2 &H^3 &H^4 &H^5 \\
\z & 0 & \z^{s-1}+\z/d & \z^{s-1}+\sum_i (\z/m_i)^{2g(D_i)} & \z/d & \z.
\end{array}
$$
\item Here $d$ is the largest natural number such that
\begin{enumerate}
\item $M(\Delta)\cdot c_1(Y/S)\in \weil(S)$ is divisible by $d=d_w$, or
equivalently
\item $M(S,\Delta)\cdot c_1(Y/S)\in \pic(S)$ is divisible by $d=d_p$.
\end{enumerate}
\end{enumerate}
\end{thm}

Proof. 
The cohomology groups $H^i(L,\z)$ are computed by a 
Leray  spectral sequence whose $E_2$ term is
$$
E^{i,j}_2=H^i(S, R^jf_*\z_{L})\Rightarrow H^{i+j}(L,\z).
$$
Every fiber  of $f$ is $S^1$, so 
$R^2f_*\z_L=0$ and
the only
interesting
higher direct image is $R^1f_*\z_{L}$.
We start by computing it.

\begin{lem} \label{H2tors=branch}
Let $L^5\to (S,\Delta=\sum_i (1-\tfrac1{m_i})D_i)$
 be a Seifert bundle over a  surface.
Assume that $H^{orb}_1(S,\Delta)=0$. Then
\begin{enumerate}
\item $H^1(S,R^1f_*\z_L)$ is torsion.
\item $H^2(S,R^1f_*\z_L)\cong \z^s+\sum_i (\z/m_i)^{2g(D_i)}$
where $s=\rank H^2(S,\q)$,
\item  $H^3(S,R^1f_*\z_L)=0$.
\item  $H^4(S,R^1f_*\z_L)=\z$.
\end{enumerate}
\end{lem}

Proof.  By \cite[48]{k-seif} there is an exact sequence
$$
0\to R^1f_*\z_L\stackrel{\tau}{\to} \z_S \to Q\to 0
\eqno{(\ref{H2tors=branch}.5)}
$$
where $Q$ in turn sits in another exact sequence
$$
0\to \sum_i \z_{P_j}/n_j \to Q \to \sum_i \z_{D_i}/m_i\to 0,
$$
where $P_j\in S$ are the singular points.
Thus  $H^0(S,Q)$ is torsion and
$H^i(S,Q)=\sum_i H^i(D_i,\z/m_i)$ for $i\geq 1$.
Putting these into the long cohomology sequence of (\ref{H2tors=branch}.5)
 we get 
(1) and (4) right away.

$H^1(S,\z)=H^1(S^0,\z)=0$ since  $H^{orb}_1(S,\Delta)=0$, and by 
(\ref{rtl.surf.coh.prop}.1) this implies that
$H^3(S,\z)=0$. 
The remaining sequence is
$$
\begin{array}{rl}
0\to \sum_i H^1(D_i,\z/m_i)&\to H^2(S,R^1f_*\z_L)
\to H^2(S,\z)\\
&\to \sum_i H^2(D_i,\z/m_i)\to H^3(S,R^1f_*\z_L)\to 0.
\end{array}
$$
$H^2(S,\z)\to \sum_i H^2(D_i,\z/m_i)$ is surjective by 
(\ref{ratsurf.orbif.simpconn.prop}),
which gives the rest.\qed
\medskip

 In  the Leray spectral sequence
$H^i(S,R^jf_*\z_L)\Rightarrow H^{i+j}(L,\z)$ the $E_2$ term is
$$
\begin{array}{lccll}
\z \quad & (\mbox{torsion}) & \quad\z^{s}+\sum_i (\z/m_i)^{2g(D_i)}\quad & 
0\quad & \z\\
\z & 0 & \z^{s} & 0 & \z.
\end{array}
$$
The spectral sequence degenerates at $E^3$ and we have only
two nontrivial differentials
$$
\delta_0: E^{0,1}_2\to E^{2,0}_2 \qtq{and} \delta_2:E^{2,1}_2\to E^{4,0}_2.
$$
By \cite[44]{k-seif}, the image of
$\delta_0:\z\to \z^s$
  is generated by $M(S,\Delta)\cdot c_1(Y/S)$.
It is nonzero since $H_1(L,\q)=0$ by assumption.
Thus $\delta_2$ is also nonzero, either by multiplicativity or by noting
that  otherwise Poincar\'e duality would fail on $L$.

Hence the $E_3$ term is
$$
\begin{array}{lcclc}
0  \quad & (\mbox{torsion}) & \quad\z^{s-1}+\sum_i (\z/m_i)^{2g(D_i)} 
\quad & 0 
 \quad & \z\\
\z & 0 & \z^{s-1}+\z/d_p & 0 & (\mbox{torsion}),
\end{array}
$$
where $d_p$ is as defined in (\ref{simpconn.over.ratsurf.thm}.2.b).
The torsion in $E^{i,0}_3$ injects into the torsion
in $H^{i}(L,\z)$.

In the notation of (\ref{fiber.trison.coeff}),
$H_1(L,\z)$ is generated by $k$ and its order  divides
$d_w$ as defined in (\ref{simpconn.over.ratsurf.thm}.2.a).
(We need to apply (\ref{fiber.trison.coeff}) to $S^0$ instead of $S$,
and $H^2(S^0,\z)\cong H_2(S,\z)$. 
Since $ H_2(S,\z)/\weil(S)$ is torsion free by (\ref{pic.weil.prop}),
divisibility in $ H_2(S,\z)$ is the same as divisibility in $\weil(S)$.)

Thus $d_p$ divides $d_w$ and 
 we need to prove that $d_p=d_w$.

$\weil(S)/\pic(S)$ is $M(S)$-torsion and
$M(S,\Delta)\cdot c_1(Y/S)= K\cdot M(\Delta)\cdot c_1(Y/S)$
for some $K|M(S)$ by (\ref{local.trosion.measures}). 
Thus divisibility by  a number relatively prime to
$M(S)$ is unchanged when we go from
$\weil(S)$ to $\pic(S)$.
Therefore we only need to prove that
 $M(\Delta)\cdot c_1(Y/S)\in \weil(S)$ is not divisible by any prime divisor
of $M(S)$.

Assuming the converse, let $p$ be such a prime
and write $M(S)=p^ap'$ where $p$ does not divide $p'$.
Assume that  $M(\Delta)\cdot c_1(Y/S)=p\cdot A$ for some $A
\in \weil(S)$.
If $p\not\vert M(\Delta)$ then $M(\Delta)\cdot c_1(Y/S)$
has order divisible by $p^a$ in some $\weil(x,X)$ by
(\ref{smoothness.crit}) and so it can not be divisible by $p$.
If $p\vert M(\Delta)$ then
$M(\Delta)\cdot c_1(Y/S)=p\cdot A$ can be rearranged to be
$$
\textstyle{\sum_{i: p\vert m_i}}
 a_i D_i=p\cdot(\mbox{integral Weil divisor}),
$$
and not all the $a_i$ are divisible by $p$.
By (\ref{p-covers.cor}) this would imply
$H^{orb}_1(S,\Delta)\neq 0$, a contradiction.
\qed

\section{Log Del Pezzo surfaces with nonrational boundary}

We see from (\ref{simpconn.over.ratsurf.thm})
that torsion in $H_2(L,\z)$, or equivalently,
torsion in $H^3(L,\z)$ is connected with
higher genus curves in the branch divisor $\Delta$
of $(S,\Delta=\sum (1-\frac1{m_i})D_i)$.

The aim of this section is to study the cases when $\Delta$
has an irreducible component of genus at least 1.
This turns out to be a  strong restriction, especially if
the $m_i$ are not very small.
A result of this type is a special case of a general principle.

\begin{say}[Ascending chain conditions] Assume that
$(S,\Delta=\sum_{i=0}^n a_iD_i)$ is a pair where we only assume that
$-(K_S+\Delta)$ is nef but we consider the case when
$a_0=1$ and $D_0$ has geometric genus $\geq 1$. The 
adjunction formula  (with a little extra work for the singularities)
says that
$$
\deg K_{D_0}\leq D_0\cdot(K_S+D_0)\leq  D_0\cdot (-\sum_{i=1}^n a_iD_i)\leq 0.
$$
Thus $D_0$ is elliptic and we also get that
$\Delta=D_0$.

The ascending chain condition principle predicts that all this
should also work if $a_0$ is close to $1$, and
  \cite[sec.5]{k-logsurf} implies  this
 for $a_0\geq 41/42$. This  is quite surprising at first sight since 
we work with singular surfaces, and the various intersection numbers
are only rational numbers.

Here we are in a rather special situation, and we get better
bounds with careful case analysis.
\end{say}

\begin{defn}\label{lDP.defn}
A {\it log Del Pezzo} surface is a pair
 $(S,\Delta)$ where 
\begin{enumerate}
\item $S$ is a normal, projective surface,
\item  $\Delta:=\sum a_iD_i$ is a linear combination
of distinct irreducible divisors with $0\leq a_i\leq 1$, and
\item $-(K_S+\Delta)$ is ample. (That is, a suitable multiple of it is
an ample Cartier divisor.)
\end{enumerate}

Thus a 2--dimensional log Fano orbifold is a
log Del Pezzo surface.

The $\q$-divisor $\Delta$ is called the {\it boundary}.
 In the orbifold  case the boundary is of the form
$\Delta=\sum (1-\tfrac1{m_i})D_i$, where the $m_i$ are natural numbers.
Such a boundary is sometimes called {\it standard}.
\end{defn}

Ideally one would like to have a classification of
all log Del Pezzo surfaces, but this seems quite out of reach.
Some kind of rough structure theorems are given in
\cite{keel-mck, shokurov}.

\begin{prop} \label{dp.genus.estimate.prop}
Let $(S,\Delta=\sum_{i=0}^n a_iD_i)$ be a log Del Pezzo surface.
 Assume that 
$D_0$ has geometric genus $\geq 1$ and $a_0\geq 1/2$.
Then $D_i$ is rational for $i\geq 1$ and 
\begin{enumerate}
\item  $g(D_0)=1$ if $a_0\geq 5/6$,
\item  $g(D_0)\leq 2$ if $a_0 \geq 3/4$,
\item  $g(D_0)\leq 4$ if $a_0\geq 2/3$.
\end{enumerate}
\end{prop}

The proof of this is an easy application of the
minimal model program for rational surfaces. In general, it is
 difficult to run the minimal model program backwards, and
it seems complicated to get a full classification, except when 
$a_0$ is close to one.

\begin{prop}\label{get.DV.DP.lem}
 Notation and assumptions as in (\ref{dp.genus.estimate.prop}).
Assume in addition that
$a_0\geq 11/12$ and $a_i\geq 1/2$ for every $i$.
Then $S$ is a Del Pezzo surface with Du Val singularities
and $\Delta=a_0D_0$.
The Picard number of $S$ is at most 8.
\end{prop}

 Del Pezzo surfaces of Picard number 1 and  with Du Val singularities
 are classified in \cite{furushima}. Adding $\p^2$ and the quadric cone,
one gets 29 types.  (Some types correspond to more than 1 surface, I
do not know how many there are up to isomorphism.)
By the results of \cite{mi-zh1}, or by checking
the conditions of (\ref{ratsurf.rho=1simply.conn.cor}),
 we obtain the following.

\begin{prop}\label{simpconn.dv.dp.prop}
Let $(S,(1-\tfrac1{m})D)$ be a log Del Pezzo surface of Picard number 1, 
with $m\geq 12$ and $g(D)\geq 1$. If $H_1(S^0,\z)=0$ then
$S$ is one of the following:
\begin{enumerate}
\item $\p^2$. Here $\pic(S)=\weil(S)$ 
and it is generated by the class of a line.
\item $Q$, the quadric cone. Here $\pic(S)=2\cdot \weil(S)$,
 the latter generated
by the lines through the vertex of the cone.
\item  $\p^2(1,2,3)$, the weighted projective plane with
weights $1,2,3$. Here $\pic(S)=6\cdot \weil(S)$, the latter generated
by the line $(x=0)$ where $x$ is the weight 1 coordinate.
\item $S_5$, a degree 5 Del Pezzo surface with a single point of index 5.
This is obtained by blowing up a flex of a smooth cubic 4 times
and then contracting the tangent line and the first 3 exceptional curves.
If $G$ is the equation of the cubic and $L$ the equation of the flex tangent
then $G/L^3$ lifts to a rational function on $S_5$ which has a 5 fold
pole along the unique line of the surface which pases through the singular
 point. The line generates $\weil(S_5)$ and $\pic(S_5)=5\cdot \weil(S_5)$.
\end{enumerate}
In all of these cases, $D$ is a smooth member of $-K_S$ and it
has genus 1.\qed
\end{prop}

The following two examples show that for $m\leq 11$ 
there are other cases as well. 

\begin{exmp}\label{P123.to.11.exmp} On $\p(1,2,3)$ consider the
divisor $\Delta=\tfrac12(x=0)+\tfrac{10}{11}(x^6+y^3+z^2=0)$.
The surface $S$ is as in
(\ref{simpconn.dv.dp.prop}) but $\Delta$ is different.  
Also, the coefficient $\tfrac{10}{11}$ can be replaced with anything smaller.
\end{exmp}

\begin{exmp} \label{E8.to.10.exmp} Let $S$ be a degree 1 Del Pezzo surface
where $|-K_S|$ has a member $D_1$ of type $II^*$ (also denoted by $\tilde E_8$)
 on Kodaira's list
(cf.\ \cite[V.7]{bhpv}). This has a unique point $P$ of multiplicity
11. Blow it up $S'\to S$ and contract the birational transform $D'_1$ of $D_1$
to get $S'\to S^*$, a log Del Pezzo surface  with Picard number 1.

Let $D_0$ be a smooth elliptic member of $|-K_S|$ and
$D^*_0$ its birational transform on $S^*$. Then
$$
\tfrac{10}{11}D_0+\tfrac1{11}D_1\equiv -K_{S}
\qtq{and}\tfrac{10}{11}D'_0+\left(1-\tfrac1{11}\right)D'_1\equiv -K_{S'}.
$$
Thus
$$
-K_{S^*}\equiv \tfrac{10}{11}D^*_0\qtq{and}
-\left(K_{S^*}+\left(1-\tfrac1{10}\right)D^*_0\right)\qtq{is ample.}
$$
This gives  an example  $(S^*, (1-\tfrac{1}{10})D^*_0)$  which is not
part of the main series.
\end{exmp}

\begin{say}[Proof of (\ref{dp.genus.estimate.prop})]
\label{pf.of.dp.genus.estimate.prop}
 Let $g:S'\to S$ be the minimal resolution of $S$.
$K_{S'}\equiv g^*K_S-(\mbox{effective divisor})$, thus
$-K_{S'}\equiv a_0D'_0+(\mbox{big effective divisor})$
where $D'_0\subset S'$ is the birational transform of $D_0$.
Since $S'$ is rational, there is a morphism
 $h:S'\to S^m$ where $S^m$ is either $\p^2$ or a minimal ruled surface.
Only rational curves are contracted by $h$,
so $-K_{S^m}\equiv a_0D_0^m+(\mbox{big effective divisor})$
where $D_0^m=h(D_0')$.

If $S^m\cong \p^2$ and  $\deg D_0^m\geq 4$, then
$4a_0\leq a_0\deg  D_0^m <3$
implies that $a_0<3/4$.
Thus if  $a_0\geq 3/4$ then $D_0^m\subset \p^2$ is a smooth degree 3 curve.

Similarly, if $S^m=\p^1\times \p^1$, the smooth quadric,
we get that $D_0^m$ is a smooth elliptic curve cut out by
another quadric if $a_0\geq 2/3$.

The remaining  possibility is that
 $S^m\cong \f_n$, the minimal ruled surface with a section $E$ of
selfintersection $-n$ and $n\geq 2$.  Let $F$ denote a fiber of 
the projection to $\p^1$.
Then $-K_{\f_n}=2E+(n+2)F$. If $|aE+bF|$ has a nonrational member
then $a\geq 2$ and $b\geq na$ and so in our case
$$
2E+(n+2)F-a_0(2E+2nF)=(2-2a_0)E+ (n+2-2na_0)F
$$
 is effective and big. Thus  $n+2>2na_0$ and 
for $a_0\geq 2/3$  this gives $n\leq 5$.
We also get that $a_0\leq 5/6$ for $n\geq 3$.

If $a\geq 3$ then the coefficient of $E$ is
$\leq 2-3a_0$, thus again $a_0<2/3$. Hence we  need to enumerate all cases
when $n=2,3,4,5$ and $C\in |2E+(2n+c)F|$ for some $c\geq 0$.
As before, we get that
$(2-2a_0)E+ (n+2-(2n+c)a_0)F$
 is effective and big which implies that $a_0< (n+2)/(2n+c)$.

If $n=2$ and $c=0$ then $D_0$ is elliptic. In all other cases we
get that $a_0<5/6$ and $a_0\geq 4/5$ only if $n=3,c=0$.
Furthermore, $a_0\geq 3/4$ in one additional case only,
when  $n=2, c=1$.
\qed
\end{say}

\begin{exmp}\label{over.ruled.surf.exmp}
 Let $F_n\subset \p^{n+1}$ be the cone over
the rational normal curve of degree $n$ in $\p^n$.
It can be also realized as $F_n=\p(1,1,n)$.
$\weil(F_n)$ is generated by the lines $L$.
$K_{F_n}\sim -(n+2)L$ and $\pic(F_n)$ is generated by
the hyperplane class which is $nL$.

Let $C\subset F_n$ be a smooth intersection of $F_n$ with a  quadric.
Thus $C\in |2nL|$ and 
$g(C)=n-1$. The surface  $(F_n, (1-\tfrac1{m})C)$
is log Del Pezzo and $g(C)>1$ in 
 the following cases.
\begin{enumerate}
\item  $n=3$  and $m\leq 5$ with $g(C)=2$
\item $n=4$ and $m\leq 3$ with $g(C)=3$.
\item $n=5$ and $m\leq 3$ with $g(C)=4$.
\item $n\geq 3$ and $m=2$ with $g(C)=n-1$
\end{enumerate}
\end{exmp}

\begin{exmp}  On $\p(1,2,5)$ consider
a general member $C\in |\o(10))|$. $C$ is smooth, has genus 2 and
the pair $( \p(1,2,5),\tfrac34 C)$ is  log Del Pezzo. 
\end{exmp}

\begin{prop}\label{cubic-cone.torsion.prop}
 Let $(S, (1-\tfrac1{m})D_0+\sum a_iD_i)$
 be a log Del Pezzo surface.
Assume that $g(D_0)\geq 2$, $m\geq 5$ and $a_i\geq \frac12$.
 Then $S\subset \p^4$ is the cone over
the rational cubic and $D_0$ is the intersection of $S$ with a quadric.

There is a unique such pair  $(S, (1-\tfrac15)D_0)$
for every genus 2 curve $D_0$.
\end{prop}

Proof. From the case analysis in (\ref{pf.of.dp.genus.estimate.prop})
we see that the minimal resolution of $S$ dominates 
$\f_3$ and does not dominate any other minimal rational surface.
Any one point blow up of $\f_3$ dominates either $\f_2$ or
$\f_4$, thus the minimal resolution of $S$ is $\f_3$.

The first possibility is that $S$ is the cone over
the rational cubic and then 
$D$ is the intersection of $S$ with a quadric by the cases analysis.

The second possibility is that $S=\f_3$.
Using the notation of (\ref{pf.of.dp.genus.estimate.prop})
we get that $(1-\tfrac1{m})D\leq -K_S=2E+5F$. Since $m\geq 5$,
this implies that $D\leq 2E+6F$. The condition $g(D)\geq 2$
implies that in fact $D= 2E+6F$. Then
$-(K_S+(1-\tfrac1{m})D)\leq \tfrac25E+\tfrac15F$ and
there is no room for any other boundary curves. Finally
$-(K_S+(1-\tfrac15)D)$ has negative intersection number with
$E$, a contradiction.

Let $D$ be any genus 2 curve and $f:D\to \p^1$ the canonical map.
Then $f_*\o_D=\o_{\p^1}+\o_{\p^1}(-3)$. This provides the unique embedding of
$D$ into $\f_3$, the projectivization of $\o_{\p^1}+\o_{\p^1}(-3)$.
\qed

\begin{say}[Proof of (\ref{get.DV.DP.lem})]
 As we noted during the proof of (\ref{dp.genus.estimate.prop}),
$a_0\geq 5/6$ implies that  $S^m=\p^2, \p^1\times \p^1 $ or $\f_2$.

$S'$ is obtained from $S^m$ by repeatedly blowing up points
$S'=S_n\to S_{n-1}\to\cdots \to S^m$.

Let us look at the first blow up $h:S_j\to S_{j-1}$
where we blow up a point $p$ not on the birational transform of $D_0$.
Let $E\subset S_j$ denote the exceptional curve of the blow up.
On $S_j$ we can write 
$$
-K_{S_j}\equiv a_0D_0^j+\alpha E+\Delta_j^*,
$$
where  $\alpha\geq 0$, 
$\Delta_j^*$ is effective and its support does not contain $E$.
By our assumption, $E$ is disjoint from $D_0^j$.
This implies that 
$$
(E\cdot \Delta_j^*)\geq -(K_S\cdot E)-\alpha (E\cdot E)=1+\alpha\geq 1,
$$
 which in turn gives
that
$$
-K_{S_{j-1}}\equiv a_0D_0^{j-1}+h_*(\Delta_j^*)
\qtq{and} \mult_ph_*(\Delta_j^*)\geq 1.
$$
$1/(1-a_0)h_*(\Delta_j^*)$ is numerically equivalent to
$-K_{S_{j-1}}$ and it has multiplicity at least 
$1/(1-a_0)\geq 12$ at $p$.
This is impossible by (\ref{weakDP.maxmult.lem}).

 Thus 
as we go from $S^m$ to $S'$,
we can blow up at most 8 points on $\p^2$
(or at most 7 points on $\p^1\times \p^1$ or $\f_2$)
and all the blow ups occur on the birational transform of $D_0$.
 Hence
$S$ is a  Del Pezzo surface with Du Val singularities and its 
Picard nunber  is at most 9.

Finally, if we had any other curve $D_i$ in $\Delta$, then
$D_i$ would occur with coefficient at least 1/2 in
$(1-\epsilon-a_0)K_S$, so we would obtain a curve
with coefficient bigger than 6 in an effective
 divisor numerically equaivalent to
$-K_S$.  This is again impossible by (\ref{weakDP.maxmult.lem}).
\qed

\begin{lem}\label{weakDP.maxmult.lem}
 Let $Y$ be one of the surfaces $\p^2, \p^1\times \p^1 $ or $\f_2$,
 $C\subset Y$ a smooth cubic. Let $f:Z\to Y$ be obtained by repeatedly
blowing up points of $C$. Let $D\subset Z$ be a divisor numerically
equivalent to $-K_Z$. Then
\begin{enumerate}
\item $\mult_zD\leq 11$ for every point $z\in Z$.
\item $\mult_AD\leq 6$ for every curve $A\subset  Z$.
\end{enumerate}
\end{lem}

Proof.  Let $Z=Y_n\to Y_{n-1}\to \cdots \to Y_j=Y$ be the sequence of blow ups
where we set $Y=Y_0$ if $Y=\p^2$ and $Y=Y_1$ if $Y=\p^1\times \p^1, \f_2$.

The push forward of $D$ to any of the $Y_i$ is again numerically
equivalent to $-K_{Y_i}$ and the multiplicity does not change
unless $z$ or $A$ is contained in the exceptional set. Thus by
induction on $n$ it is sufficient to consider the case when
$g:Z\to Z':=Y_{n-1}$ is the blow up of a point $p$ and 
 $z$ or $A$ is contained in the exceptional curve $E$ of $g$.
Set $D':=g_*(D)$, then $D=g^*(D')+(\mult_pD'-1)E$. 

If $n\geq 10$ then $(-K_Z)^2=9-n<0$ and so the only divisor
numerically equivalent to $-K_Z$ is $C$ itself which is smooth.

Let $C'\subset Z'$ denote the birational transform of $C$. Then
$\mult_pD'\leq (C'\cdot D')=9-(n-1)$, hence if $n\geq 4$ then
$\mult_pD'\leq 6$ and
so $\mult_ED\leq 5$ and $\mult_zD\leq 11$ for every $z\in E$.

We are left with the cases when $n\leq 3$, that is, we blow up at most
$3$ points. In these cases $Z$ is a toric variety. By degenerating
$D$ using the torus action and noting that multiplicty is upper semi
continuous in degenerations, we are reduced to the case when
$D$ is fixed by the torus. These cases are easy to enumerate by hand.
\qed

\end{say}

\section{Log Del Pezzo surfaces with $H_1^{orb}(S)=0$}

The aim of this section is to classify
Del Pezzo surfaces with cyclic Du Val singularities
satisfying  $H_1^{orb}(S)=0$.

There are two ways to proceed.

\begin{say}[Traditional approach]
The minimal resolution of a  Del Pezzo surface $S$ with  Du Val singularities
 is obtained from $\p^2$
and a smooth cubic curve $C\subset \p^2$ by
blowing up  $m\leq 8$ points on the cubic.

In the smooth case, the $m$ points are different and we get
only one family for each $m$. In the singular case,
the deformation type of the resulting surface is determined
by the following data:
\begin{enumerate}
\item Any number of the $m$ points may coincide, and we have to
mark all coinciding point pairs.
\item We have to mark  triplets of  points that are on a line.
(Keep in mind that if 3 points coincide, they are also on a line
if this point is  an inflection point of the cubic.)
\item We have to mark  sextuplets of  points that are on a conic.
(Again, this can happen with 6 points coinciding at a
6--torsion point.)
\end{enumerate}
It is clear that the number of cases is finite, but 
it would be quite tedious to get a complete list.

The main problem with this approach is that I do not see any efficient way
to compute  $H_1^{orb}(S)$. This seems to be nonzero for the
majority of all cases.
\end{say}

\begin{say}[Minimal model approach] Here we start with
a Del Pezzo surface $S$ with  Du Val singularities
and run the minimal model program
to get $g:S\to S^m$.
We clasify $S$ using the following observations.
\begin{enumerate}
\item  $H_1^{orb}(S)$  and  $\pi_1^{orb}(S)$  are
unchanged during the minimal model program; see
(\ref{H^1.orb.unch.mmp}).
\item $S^m$ is either $\p^1\times \p^1$ or one of the 4 surfaces
listed in (\ref{simpconn.dv.dp.prop}); see (\ref{5.min.mods}).
\item For  a fixed $S^m$, the deformation types of all possible surfaces
$S\to S^m$ are classified by sequences of natural numbers
$m_1\leq \cdots \leq m_k$ such that $\sum m_i<(K_{S^m})^2$;
see (\ref{bu.notation.dp}).
\item Many of the surfaces have representations over
different surfaces $S^m$, but these are not too hard to understand;
see (\ref{simpconn.cyc.DV.thm}).
\end{enumerate}
\end{say}

\begin{say}[Contractions with Du Val singularities]\label{H^1.orb.unch.mmp}
Let $T$ be any projective surface with Du Val singularities
and $f:T\to U$ a birational contraction 
with connected exceptional curve $E\subset T$
such that $-K_T$ is $f$-ample. These are easy to  classify,
see \cite[3.3]{keel-mck}. One gets that  $u:=f(E)\in U$ is a smooth point
and $f$ can be described as follows.

Let $u\in C\subset U$ be a smooth curve germ. Fix a number
$m$ and blow up $u\in C$ repeatedly $m$-times. We get $m$ exceptional
curves. One of these is a $-1$-curve, the other $m-1$ form a chain of
$-2$-curves. This chain can be contracted to a point of type
$A_m$.   (It can be described by local equations as $xy-z^{m+1}=0$
or as a quotient of  $\c^2$ by the $\z/(m+1)$-action
$(x,y)\mapsto (\epsilon x, \epsilon^{-1} y)$.) The resulting 
surface is $T$. Thus we see that $E\subset T$ is orbismooth,
$T$ has 
 a unique singular point  $t\in E$ and a neighborhood of
$E$ minus $t$ is homotopic to a lense space $S^3/(\z/m)$
with a disc attached (corresponding to $E\setminus \{t\}$),
killing the fundamental group. From this, or using
\cite[3.3]{keel-mck} we conclude that
$H_1^{orb}(T)=H_1^{orb}(U)$ and
$\pi_1^{orb}(T)=\pi_1^{orb}(U)$.

The deformation type of $T$ is determined by the deformation type of
$U$ and the number $m$.
\end{say}

\begin{defn}\label{bu.notation.dp}
 Let $S$ be a Del Pezzo surface with Du Val singularities
and $m_1\geq \cdots \geq m_k\geq 1$ integers. We denote by
$B_{m_1,\dots,m_k}S$ any surface obatined as follows.

Pick any smooth elliptic curve $C\in |-K_S|$ and $p_1,\dots,p_k$ distinct
points on $C$. Then perform a blow up type $m_i$ at $p_i$.

All such surfaces form one deformation type.
Furthermore, $S$ and the deformation type
determine the numbers $m_i$. Indeed, 
the numbers $m_i\geq 2$ can be read off from
the singularities    and the Picard number
 determines $k$.

The canonical class of $B_{m_1,\dots,m_k}S$
is nef and big iff $\sum m_i<(K_S^2)$. If this holds then
 $B_{m_1,\dots,m_k}S$ is a Del Pezzo surface for general
choice of the points $p_i$.
\end{defn}

\begin{say}[The minimal models $S^m$]\label{5.min.mods}
If we start with a Del Pezzo surface  with (cyclic) 
Du Val singularities $S$, then all steps
of the minimal model program yield Del Pezzo surfaces with (cyclic) 
Du Val singularities.
The program eventually stops, and we end up with
$S\to S^m$ where $S^m$ does not have birational contractions.
This can happen in two ways
\begin{enumerate}
\item $S^m$ has Picard number 1. If, in addition,
$H_1^{orb}(S)=0$ and $S$ has cyclic 
Du Val singularities then $H_1^{orb}(S^m)=0$ and 
by (\ref{simpconn.dv.dp.prop})  $S^m$ is one of the 4 surfaces 
$$
\p(1,2,3), Q, \p^2 \qtq{or} S_5.
$$
\item $S^m$ has Picard number 2 and it has two different
conic bundle sutructures. For $S^m$ smooth, this happens
only for $\p^1\times \p^1$.
A list of the  singular cases is given in \cite{mi-zh2}.
They are either not simply connected or have noncyclic
singularities.
\end{enumerate}
\end{say}

We can now state the main classification theorem of the section:

\begin{thm}\label{simpconn.cyc.DV.thm}  There are 93 deformation types of
Del Pezzo surfaces with cyclic Du Val singularities
satisfying  $H_1^{orb}(S)=0$.
These are
\begin{enumerate}
\item $B_{m_1,\dots,m_k}\p(1,2,3)$ for $\sum m_i<6$ and $m_i\geq 2$,
\item $B_{m_1,\dots,m_k}Q$ for $\sum m_i<8$ and $m_i\geq 2$,
\item $B_{m_1,\dots,m_k}\p^2$ for $\sum m_i<9$ and $m_i\geq 2$,
\item $B_{m_1,\dots,m_k}\p^1\times \p^1$ for $\sum m_i<8$,
\item  $S_5, B_3S_5, B_4 S_5$ and $B_1\p^2$.
\end{enumerate}
\noindent All these satisfy  $\pi_1^{orb}(S)=0$.
\end{thm}

Proof. We know that all of the surfaces are of the form
 $B_{m_1,\dots,m_k}S$ where $S$ is one of the 5 surfaces
listed in (\ref{5.min.mods}). First we prove that if
$m_k=1$ and $S\neq \p^1\times \p^1$ then the
resulting surface is somewhere else on the list.
This follows from the easy isomorphisms:
$$
B_1\p(1,2,3)\cong B_3Q, B_1Q\cong B_2\p^2, B_{m,1}\p^2\cong 
B_m\p^1\times \p^1, B_1S_5\cong B_4\p^2.
$$
This implies that  every 
Del Pezzo surface with cyclic Du Val singularities
satisfying  $H_1^{orb}(S)=0$ is either on the list or
it is $B_{m_1,\dots,m_k}S_5$ for $\sum m_i<5$ and $m_i\geq 2$.
There are only 5 such cases. for two of these, $m_k=2$.
These are eliminated by the harder isomorphism
$$
B_2S_5\cong B_4 Q.
$$
The easiest way to see this is to notice that these blow ups give a 
cubic surface,  whose explicit model is
described in \cite[Species XIV, p.79]{henderson}.
In his notation, contracting the line $x=z=0$ gives $Q$ and
contracting the line $y=w=0$ gives $S_5$.

It remains to see that all these are different.
Amost all cases are decided by looking at the
number of singular points minus the Picard number.
For the surfaces in (\ref{simpconn.cyc.DV.thm}.1)
we get 1,
for those in (\ref{simpconn.cyc.DV.thm}.2)
we get 0, (\ref{simpconn.cyc.DV.thm}.3) gives -1
and  (\ref{simpconn.cyc.DV.thm}.4) gives $\leq -2$.

It remains to observe that the 4 surfaces in
(\ref{simpconn.cyc.DV.thm}.5) are nowhere else on the list.\qed

\begin{rem} The results of \cite{mi-zh2}, coupled with
(\ref{simpconn.cyc.DV.thm}) give a complete list of all
 deformation types of
Del Pezzo surfaces with arbitrary Du Val singularities
satisfying  $H_1^{orb}(S)=0$.

By \cite[3.10]{keel-mck} and \cite{mi-zh2}, there are
5 more such surfaces with Picard number 1, with singularities
$D_5,E_6,E_7,E_8,E_8$ and 3 more with Picard number 2,
with 2 conic bundle structures. These have  singularities
$D_4,D_6,D_7$.
The blow ups considered in
(\ref{bu.notation.dp}) introduce only cyclic quotient
singularities, hence there are no ismorphisms between
the blow ups, except possibly for the two surfaces of type $E_8$.
These, however, have $K^2=1$, so we can not blow up at all.
\end{rem}

\section{Einstein metrics on Seifert bundles}

The main impediment to apply
(\ref{exst.of.E.metric}) is the current shortage of
existence results for positive Ricci curvature 
K\"ahler--Einstein metrics on orbifolds.

We use the following sufficient algebro--geometric condition.
There is every reason to expect that it is very far from being optimal,
but it does provide a large selection of good exmaples.

In this paper we use (\ref{nadel.thm}) only for surfaces.
The concept {\it klt} is defined in (\ref{klt.etc.defn}).

\begin{thm}\cite{nadel, dk}
\label{nadel.thm} Let $(X,\Delta)$ be an $n$-dimensional compact
orbifold such that $-(K_X+\Delta)$ is ample. Assume that
there is an $\epsilon>0$ such that 
$$
(X,\Delta+\tfrac{n+\epsilon}{n+1}D) \qtq{is klt}
$$
for every effective $\q$-divisor $D\equiv -(K_X+\Delta)$. Then
$(X,\Delta)$ has an orbifold  K\"ahler--Einstein metric.\qed
\end{thm}

\begin{defn}\label{klt.etc.defn}
 Let $X$ be a normal complex space  and $D$ a $\q$-divisor on $X$.
Assume that $mK_X,mD$ are both Cartier for some $m>0$
(this is automatic if $X$ is an orbifold).
Let $g:Y\to X$ be any proper birational morphism, $Y$ smooth.
 Then there is a unique
$\q$-divisor $D_Y=\sum e_iE_i$ on $Y$ such that 
$$
K_Y+D_Y\equiv g^*(K_X+D)\qtq{and}  g_*D_Y=D.
$$
We say that $(X,D)$ is {\it klt}
 ( resp.\ {\it log canonical}) if
$e_i< 1$ (resp.\ $e_i\leq 1$) for every $g$
and for every $i$.
We say that $(X,D)$ is
{\it canonical}  if
$e_i\leq 0$  for every $g$
and for every $i$ such that $E_i$ is $g$-exceptional
\end{defn}

It is quite hard to check using the above 
 definition if a pair $(X,D)$ is klt or not.
 For surfaces, there are reasonably sharp
multiplicity conditions  which ensure that a given pair $(X,D)$ 
is klt. These conditions are not  necessary, but
they seem to apply in most cases of interest to us.

\begin{say}[How to check if $(X,D)$ is klt or not?]{\ }
\label{how.to.check}

Let $X$ be a surface with quotient singularities.
Let the singular points be $P_i\in X$ and we write these
locally  as 
$$
p_i:B^2 \to B^2/G_i
$$
where $B^2$ is the unit 2--ball  and $G_i\subset GL(2,\c)$
 a finite subgroup. We may assume that the origin is an
 isolated fixed point of every
nonidentity element of $G_i$ (cf. \cite{briesk}). 
Let $D$ be an effective
$\q$-divisor on $X$.
  Then  $(X,D)$ is klt if the following three
conditions  are satisfied.
\begin{enumerate}
\item  (Klt along curves) $D$ does not contain an
irreducible component with coefficient  $\geq 1$.
\item (Klt  at smooth points) For every smooth point $P\in X$,
\begin{enumerate}
\item either $\mult_PD\leq 1$, 
\item or $D=cC+D'$ where $C$ is a curve through $P$, smooth at $P$,
$D'$ is effective not containing $C$,
and the local intersection number $(C\cdot_P D')<1$.
\end{enumerate}
 (This follows from \cite[4.5 and 5.50]{kmbook}.)
\item (Klt  at singular points) The condition (2) is satisfied
for $(B^2, p_i^*D)$ at $P=(\mbox{origin})$.
 (This follows from \cite[5.20]{kmbook} and the previous case.)
\end{enumerate}
\end{say}

A good illustration of how to use these methods is given by
the following example. Some of its conditions seem artificial,
but they are satisfied in many cases.

\begin{lem}\label{KE.ex.conds.lem}
 Let $S\subset \p^d$ be a Del Pezzo surface 
(with quotient singularities) of degree $d$
with hyperplane section  $H$.
Let $D\subset S$ be a smooth divisor, not passing through any singular points.
Assume that
\begin{enumerate}
\item $-K_S\equiv a H$ for some  $a\in\q$,
\item $D\equiv b H$ for some  $b\in\q$,
\item every singular point is a  quotient by a group of order
at most $d$,
\item there is a line $L\subset S$ passing through all the singular points
such that $L\equiv \frac1{d}H$, 
\item $d(a-(1-\frac1{m})b)<\frac3{2+\epsilon}$
 and $bd(a-(1-\frac1{m})b)<\frac3{2+\epsilon}$.
\end{enumerate}
Then    $(S,(1-\frac1{m})C)$ 
has an orbifold  K\"ahler--Einstein metric.
\end{lem}

Proof. Let $D\subset S$ be any effective divisor
numerically equivalent to
$$
-(K_S+(1-\tfrac1{m})C)\equiv (a-(1-\tfrac1{m})b)H.
$$
We need to check the conditions of
(\ref{how.to.check}) for $(S, (1-\tfrac1{m})C+\frac{2+\epsilon}3 D)$.

Since $D\cdot H=d(a-(1-\frac1{m})b)<\frac3{2+\epsilon}$, we see that
$\frac{2+\epsilon}3 D$ can not contain any effective curve with
coefficient $\geq 1$.
For the same reason, $D$ has multiplicity $\leq 1$ at every point.

At a point on $C$, we need that
$\frac{2+\epsilon}3 D\cdot C<1$. 
This follows from $bd(a-(1-\frac1{m})b)<\frac3{2+\epsilon}$.

Finally consider a singular point, where $p:B^2\to S$ has degree $d'\leq d$.
Here 
$$
\mult p^*D\leq (p^*D\cdot p^*L)=d'(D\cdot L)\leq d'\tfrac1{d}(D\cdot H)
<\tfrac3{2+\epsilon},
$$
if $L\not\subset \supp D$. 
If $D=cL+D'$ where $L\not\subset  \supp D'$, then we can use the same
estimate for $D'$ to conclude the proof.\qed

\begin{exmp}[Existence of K\"ahler--Einstein metrics]\label{KE.ex.exmp}
Here we see what (\ref{KE.ex.conds.lem}) gives for the 
surfaces that we need for (\ref{RHS.main.thm}) and (\ref{RHS.main2.thm}.2).

Let us start with the ones in (\ref{simpconn.dv.dp.prop}).

(1) $(\p^2, (1-\frac1{m})D)$ where $D$ is a smooth cubic.
Here $d=1, a=b=3$ and we get a K\"ahler--Einstein metric for $m>6$.

(2) $(Q, (1-\frac1{m})D)$ where $Q\subset \p^3$ is a quadric cone
and  $D$ is its intersection with a quadric.
Here $d=2, a=b=2$ and we get a K\"ahler--Einstein metric for $m>5$.

(3) $(\p(1,2,3), (1-\frac1{m})D)$ with an embedding
$\p(1,2,3)\subset \p^6$ as a degree 6 surface and 
$D$ is a smooth hyperplane section.
The curve $((\mbox{weight 1 coordinate})=0)$ is the required line.
Here $d=6, a=b=1$ and we get a K\"ahler--Einstein metric for $m>4$.

(4) $(S_5, (1-\frac1{m})D)$ with an embedding
$S_5\subset \p^5$ as a degree 5 surface and 
$D$ is a smooth hyperplane section.
The needed line is constructed in (\ref{simpconn.dv.dp.prop}).
Here $d=5, a=b=1$ and we get a K\"ahler--Einstein metric for $m>3$.
\medskip

Finally let us consider some of the surfaces from (\ref{over.ruled.surf.exmp})
 and (\ref{cubic-cone.torsion.prop}).
\medskip

(5)  Let $F_n\subset \p^{n+1}$ be the cone over the degree $n$ rational normal
curve and $C\subset F_n$ a smooth intersection of $F_n$ with a quadric.
Here $d=n, a=1+\frac2{n}$ and $b=2$.  The conditions of  
(\ref{KE.ex.conds.lem}) are satisfied and 
 we  get a K\"ahler--Einstein metric
in the following cases:
$(F_3, \frac45 C), (F_3, \frac34 C),  (F_4, \frac23 C)$ and $(F_5, \frac23 C)$.
\end{exmp}

\begin{lem}\label{e-metric.dv.dp}
 Let $S$ be a Del Pezzo surface with Du Val singularities,
$C\in |-K_S|$ a smooth elliptic curve and $D$ an effective divisor on $S$
numerically equivalent to $\frac1{m}(-K_S)$. Assume that $m\geq 9$.

Then $(S,(1-\frac1{m})C+D)$ is log canonical,
thus $(S,(1-\frac1{m})C)$ has an orbifold K\"ahler--Einstein metric.
\end{lem}

Proof. Since $K_S+(1-\frac1{m})C+D$ is numericaly zero,
being log canonical is preserved under pull backs and push forwards.
Thus we can pass first to the minimal resolution, and then
go down to $S=\p^2$ or $S=\p^1\times\p^1$. These are now easy
to treat with the estimates of (\ref{how.to.check}).\qed

\section{Seifert bundles and  rational homology spheres}

With the aim of constructing new Einstein metrics on homology spheres,
 one would like to describe all Seifert bundle structures on them.

In view of (\ref{rho=1.orbH1=0.cor}), the following is a
more detailed version of (\ref{rhs.seif-strs.thm}).

\begin{thm} \label{seif.on.S5.thm}
Let $f:L\to (S,\sum (1-\frac1{m_i})D_i)$ be a 5-dimensional
Seifert bundle, $L$ smooth.
\begin{enumerate}
\item If $L$ is a rational homology sphere with $H_1(L,\z)=0$
then
\begin{enumerate}
\item $S$ has only cyclic quotient
singularities and $H_2(S,\z)\cong \weil(S)\cong \z$,
\item the $D_i$ are orbismooth  curves, intersecting
transversally,
\item the $m_i$ are relatively prime to each other and 
each $m_i$ is relatively prime to  $\deg D_i\in H_2(S,\z)\cong\z$.
\end{enumerate}
\item 
Coversely, given any $(S,\sum (1-\frac1{m_i})D_i)$ satisfying the above
3 conditions (1.a,b,c), there is a unique (up to orientation)
Seifert bundle $f:L\to (S,\sum (1-\frac1{m_i})D_i)$
such that $L$ is a rational homology sphere with $H_1(L,\z)=0$.
\item  $L$ is simply connected iff $\pi^{orb}_1(S,\sum (1-\frac1{m_i})D_i)=1$.
\item $L$ is homeomorphic to $S^5$ iff
$\pi^{orb}_1(S,\sum (1-\frac1{m_i})D_i)=1$ and the $D_i$ are all rational.
\end{enumerate}
\end{thm}

This gives a rather
  complete answer in terms of algebraic geometry.
There is, however, one missing piece.
If (\ref{simpconn.implies.rtl.conj}) is true then
we get further restrictions of $S$:

\begin{conj} \label{seif.on.S5.conj}
Let $f:L\to (S,\sum (1-\frac1{m_i})D_i)$ be a 5-dimensional
Seifert bundle, $L$ smooth.
 If $L$ is a simply connected rational homology sphere 
then $S$ is a rational surface.
\end{conj}

\begin{rem} A rich source of examples, first considered in \cite{or-wa},
comes from taking $S=\p^2$. If the  $D_i$ are lines and conics
intersecting transversally, the  $m_i$ are relatively prime 
to each other  and
odd for conics, then  $L$ is always $S^5$.

By allowing higher degree curves for the $D_i$, we get
many examples of simply connected rational homology spheres.
However, not all  rational homology spheres can be realized.
Indeed, the torsion subgroup of $H_2(L,\z)$
is computed by (\ref{simpconn.over.ratsurf.thm}).
 If $H_2(L,\z)$ is torsion, then
$S$ has Picard number one, hence any two curves on $S$ intersect.
By (\ref{loc.cyc.orb.say}.2)
 this implies that the $m_i$ are relatively prime 
to each other. Hence we obtain:
\end{rem}

\begin{cor} \label{rhs.seif-restriction}
Let $f:L\to (S,\sum (1-\frac1{m_i})D_i)$ be a 5-dimensional
Seifert bundle, $L$ smooth.
 If $L$ is a rational homology sphere with $H_1(L,\z)=0$
then 
$$
H_2(L,\z)\cong \sum_i (\z/m_i)^{2g(D_i)}
$$
where the $m_i$ are relatively prime 
to each other. \qed
\end{cor}

Thus, for instance, $H_2(L,\z)$ can not be $(\z/p)^2+(\z/p^2)^2$.

\begin{say}[Proof of (\ref{seif.on.S5.thm})]
 Let us start with $f:L\to (S,\sum (1-\frac1{m_i})D_i)$
such that  $L$ is a rational homology sphere. 
By (\ref{seif.gives.orbifold}) and (\ref{loc.cyc.orb.say}),
 $S$ has only cyclic quotient
singularities, the $D_i$ are orbismooth and they  intersect
transversally. The  rest of the
conditions follows from (\ref{rho=1.orbH1=0.cor}).

Conversely, start with a  surface $S$ with $H_2(S,\z)\cong \weil(S)\cong \z$,
 with singular points
$P_i$ with links $M_i$ with $H^2(M_i,\z)\cong \z/n_i$
 and orbismooth rational curves
$D_i$ of degre $d_i$ and natural numbers $m_i$.

The $n_i$ are pairwise relatively prime by 
(\ref{ratsurf.rho=1simply.conn.cor}); set $N=\prod n_i$.
 The $m_j$ are pairwise relatively prime by assumption; set $M=\prod m_j$.

 $H^{orb}_1(S,\sum (1-\frac1{m_i})D_i)=1$ by (\ref{rho=1.orbH1=0.cor}),
hence  by (\ref{c_1.gives.Sb.for.simpconn}),
 a Seifert bundle over $(S,\Delta)$ is uniquely determined
by its Chern class. Furthermore, by (\ref{simpconn.over.ratsurf.thm}.2),
 $H_1(L,\z)=1$
iff $c_1(L/S)=\pm\frac1{M}[\ell]$ where $[\ell]\in \weil(S)$ is the
positive generator.

If $D_j$ does not pass through $P_i$ then $n_i$ divides $d_j:=\deg D_j$
since $\weil(S)\to \weil(P_i,S)$ is surjective by 
(\ref{ratsurf.rho=1simply.conn.cor}) and its kernel is
precisely those curves whose degree is divisible by $n_i$.
Hence in this case 
$m_j$ is relatively prime to $n_i$. If  $D_j$ does pass through $P_i$
then by (\ref{loc.cyc.orb.say}),
 the  multiplicity of the fiber above $P_i$ is divisible by
$n_im_j$. Thus $m(X,\Delta)=NM$.

By (\ref{seif.on.S5.thm}.1.c) the integers $M, d_jM/m_j$ are relatively prime,
thus
$$
b_M+\sum \tfrac{d_j}{m_j}b_j= \tfrac1{M}
$$
is solvable in integers $b_M,b_j$. We can even assume that
$0\leq b_j<m_j$ and $b_j$ is necessarily relatively prime to $m_j$.
We identify $b_M\in \z$ with a Weil divisor class
$B_M\in \weil(S)\cong \z$. These data determine a
Seifert bundle $f:L\to (S,\sum (1-\frac1{m_i})D_i)$.
(The solution corresponding to $-\tfrac1{M}$ gives a Seifert
bundle which differs from this by reversing the orientation of the
cirles.)

$M\cdot c_1(L/S)$ is a generator of $\weil(S)$, and
it generates  the local class groups $\weil(s,S)$
by (\ref{ratsurf.rho=1simply.conn.cor}.3).
Thus $L$ is smooth by (\ref{smoothness.crit}).

$L$ is  a rational homology sphere by (\ref{simpconn.over.ratsurf.thm}).

 Finally,
$\pi_1(L)=1$ iff  $\pi^{orb}_1(S,\sum (1-\frac1{m_i})D_i)=1$
by (\ref{fund.gr.sequence}).

$L$ is homeomorphic to $S^5$ iff it is a simply connected
integral homology sphere. By (\ref{simpconn.over.ratsurf.thm}),
$H_2(L,\z)=0$ iff the $D_i$ are rational.
 Duality  gives  $H_3(L,\z)=0$. 
\qed
\end{say}

We are now ready to prove the main theorems of this paper.

\begin{say}[Proof of (\ref{h2.tors.thm})]

Let $f:L\to (S,\sum(1-\frac1{m_i})D_i)$ be a Seifert bundle
with $H_1(L,\z)=0$ such that 
$(S,\sum(1-\frac1{m_i})D_i)$ is a log Del Pezzo orbifold.
If $\tors H_2(L,\z)\neq 0$, then
by (\ref{simpconn.over.ratsurf.thm}.1), at least one of the curves,
say $D_0$, is nonrational. By (\ref{dp.genus.estimate.prop}),
all the others are rational and the relationships
between the genus of $D_0$ and $m_0$ given in
(\ref{dp.genus.estimate.prop}) give the
cases (\ref{h2.tors.thm}.1--4).

The existence of K\"ahler--Einstein metrics on the following
surfaces was established  in (\ref{KE.ex.exmp}.5):

$(F_3, \frac45 C)$, giving $H_2(L,\z)=(\z/5)^4$,
$(F_4, \frac23 C)$  giving $H_2(L,\z)=(\z/3)^6$, and $(F_5, \frac23 C)$
 giving $H_2(L,\z)=(\z/3)^8$.
$ (F_3, \frac34 C)$  is not orbifold simply connected, so for
$H_2(L,\z)=(\z/4)^4$ we need another example.
It is given by $(Q,\frac34 C)$ where  $Q\subset \p^3$ is 
the quadric cone and $C\subset Q$ is
a smooth degree 5 curve (thus it has to pass through the vertex).
The corresponding Seifert bundle can also be realized as the
link of the singularity $x^5+y^5+xz^2+u^4=0$.

For examples with $H_2(L,\z)=(\z/m)^2$ it is probably easiest
to use a general degree 1 smooth Del Pezzo surface $S$ and
a smooth elliptic member of $|-K_S|$. The conditions
of (\ref{how.to.check}) are very easy to check for any $m\geq 2$.

Let $\f_n$ be the minimal ruled surface with a negative section
$E\subset \f_n$ with $E^2=-n$
and fiber $F$. Take a smooth curve $C\subset |2E+(2n+3)F|$
which is transversal to $E$.
Then $g(C)=n+2$ and
$(\f_n, (1-\frac12) C+ (1-\frac1{m})E)$ is a log Del Pezzo surface
for $m>2n$.  The conditions of (\ref{how.to.check}) again present no problems
for $m\geq 7$.
These give examples with $H_2(L,\z)=(\z/2)^{2n}$ 
 for any $n\geq 2$.

Finally, an example with $H_2(L,\z)=(\z/3)^4$ 
can be obtained from  $(\f_1, \frac23 D+\frac12 E)$
where $D\in |2E+4F|$.
\qed
\end{say}

\begin{say}[Proof of (\ref{RHS.main.thm})]
\label{main.series.h2.tors}

By (\ref{seif.on.S5.thm}), the existence of 
such Seifert bundles is reduced to a question about
log Del Pezzo surfaces $(S, (1-\frac1{m})D)$ where $m\geq 12$
and $g(D)=1$.
These are classified in (\ref{simpconn.dv.dp.prop}).
We get 4 cases:
\begin{enumerate}
\item  $(\p^2, (1-\frac1{m})D)$ where $D$ is a smooth cubic.
We also need $3\not\vert m$ by (\ref{seif.on.S5.thm}.1.c).

\item  $(Q, (1-\frac1{m})D)$ where $Q\subset \p^3$ is a quadric cone
and  $D$ is its intersection with a quadric.
We also need $2\not\vert m$ by (\ref{seif.on.S5.thm}.1.c).

\item $(\p(1,2,3), (1-\frac1{m})D)$ where $D\in |\o(6)|$ is smooth.
$[D]$ is 6 times the generator in $\weil(S)$, thus
we also need $m$ to be relatively prime to 6 by (\ref{seif.on.S5.thm}.1.c).

\item $(S_5, (1-\frac1{m})D)$ with an embedding
$S_5\subset \p^5$ as a degree 5 surface and 
$D$ is a smooth hyperplane section.
$[D]$ is 5 times the generator in $\weil(S)$, thus
we also need $5\not\vert m$ by (\ref{seif.on.S5.thm}.1.c).
\end{enumerate}

The existence of K\"ahler--Einstein metrics is established in
(\ref{KE.ex.exmp}.1--4).
\qed
\end{say}

\begin{say} Let us compute $\pi_1^{orb}(S,\Delta)$
for each of the 4 cases listed in
(\ref{main.series.h2.tors}).

Let $\ell\subset S$ be the line, then
$D\sim n\cdot \ell$ where $n=3,4,6$ or $5$.
Thus $\o_S(n\cdot \ell)\cong \o_S(D)$ gives an $n$-sheeted
ramified cyclic cover $h: T\to S$ and
$K_T=h^*(K_S+(1-\frac1{n})D)$. Thus we get that
$-K_T$ is ample and 
$K_T^2=3,2,1,1$ in these 4 cases. $T$ is a smooth Del Pezzo surface
and $C:=\red h^{-1}(D)$ is a smooth elliptic curve.

If $C\subset \p^2$ is a smooth cubic then
$\pi_1(\p^2\setminus C)=\z/3$ but as soon as we blow up at least
one point on $C$, the surface will contain  a line intersecting
$C$ in one point only. Thus $\pi_1(T\setminus C)=1$ in all
4 cases and  we conclude that
$\pi_1^{orb}(\p^2\setminus D)=\z/3$,
$\pi_1^{orb}(Q\setminus D)=\z/4$,
$\pi_1^{orb}(\p(1,2,3)\setminus D)=\z/6$ and
$\pi_1^{orb}(S_5\setminus D)=\z/5$.

Our restrictions on $m$ specify that it be relatively prime to
the order of the corresponding $\pi_1^{orb}$, thus
$\pi_1^{orb}(S,\Delta)=1$ in all 4 cases.
\end{say}

\begin{say}[Proof of (\ref{RHS.main2.thm})]
\label{pf.of.8}

As before, (\ref{get.DV.DP.lem}) reduces the first part to the
 enumeration of 
all Del Pezzo surfaces $S$ with cyclic quotient Du Val singularities
such that $H_1^{orb}(S)=0$.
This is accomplished in (\ref{simpconn.cyc.DV.thm}).
K\"ahler--Einstein metrics are
constructed in (\ref{e-metric.dv.dp}).

The second part of (\ref{RHS.main2.thm})
follows from (\ref{cubic-cone.torsion.prop}) and
 (\ref{KE.ex.exmp}.5).

This completets the proof of (\ref{RHS.main2.thm}).
\qed
\end{say}

\section{Links of log terminal singularities}

Let $f:Y\to (X,\Delta)$ be a Seifert $\c^*$-bundle.
One can naturally compactify it by adding
a zero and infinity section, see \cite[14]{k-seif}.
The infinity section is contractible to a singular
point iff $c_1(Y/X)$ is ample.
As noted by \cite{pinkham, demazure}, this establishes an
equivalence between Seifert $\c^*$-bundles with $c_1(Y/X)$  ample
and singularities with a good $\c^*$-action.

The canonical class of this singularity is $\q$-Cartier
iff  the orbifold canonical class is a rational multiple of
$c_1(Y/X)$. Furthermore, the generalized adjunction formula
(cf. \cite[Sec.16]{k-etal}) implies that the singularity is 
log terminal iff  $(X,\Delta)$ is log terminal and
$-(K_X+\Delta)$ is ample.

Thus we have established an
equivalence between:
\begin{enumerate}
\item Smooth Seifert $\c^*$-bundles $f:Y\to (X,\Delta)$ such that
$-(K_X+\Delta)$ is ample and $c_1(Y/X)$ is a positive multiple of 
$-(K_X+\Delta)$.
\item Isolated log terminal singularities with a good $\c^*$-action.
\end{enumerate}

It is natural to ask, to what extent the results of this note generalize
to arbitrary log terminal singularities.

\begin{prob} Let $0\in X$ be an isolated log terminal singularity
with link $M$.
\begin{enumerate}
\item Is $\pi_1(M)$ finite?  
\item In dimension 3, do the restrictions of
(\ref{h2.tors.thm}) also apply to $M$?
\item Is there any connection between the log terminality of
$X$ and the existence of positive Ricci curvature metrics on $M$?
\item Can one obtain  Einstein metrics on links without $\c^*$-action?
\end{enumerate}
\end{prob}

Log canonical singularities with $\c^*$-action lead to
Seifert $\c^*$-bundles whose base is a Calabi--Yau
orbifold:

\begin{prop} \label{seif.over.CY.prop}
Let $f:L\to (X,\sum (1-\frac1{m_i})D_i)$ be a 
Seifert bundle, $L$ smooth.
Assume that  $H_1(L,\z)=0$ and $(X,\sum (1-\frac1{m_i})D_i)$
is a Calabi--Yau
orbifold, that is, $K_X+\sum (1-\frac1{m_i})D_i$ is numerically trivial.

Then $\sum (1-\frac1{m_i})D_i=0$ and $X$ is a Calabi--Yau
orbifold with trivial canonical class.
\end{prop}

Proof. If  $\sum (1-\frac1{m_i})D_i\neq 0$, pick a prime
$p$ dividing at least one of the $m_i$. Let us 
clear denominators in $\sum (1-\frac1{m_i})D_i\equiv -K_X$
and rewrite this as
$$
\textstyle{\sum_{i: p\vert m_i}}
 a_i D_i\sim  p\cdot(\mbox{integral Weil divisor}),
$$
where not all the $a_i$ are divisible by $p$.
Thus by (\ref{p-covers.cor}) we get that
$H_1^{orb}(X,\sum (1-\frac1{m_i})D_i)\neq 0$, a contradiction to
$H_1(L,\z)=0$.

Thus the canonical class is torsion. Thus it gives 
a torsion element of $H^2(X\setminus\sing X,\z)$ which is, however, 
torsion free since $H_1(X\setminus\sing X,\z)=0$. 
Thus the canonical class is trivial.
\qed
\medskip

In dimension 5, one can be even more precise,
as  conjectured by  Galicki:

\begin{cor}  \label{seif.over.CY.cor}
Let $f:L\to (S,\sum (1-\frac1{m_i})D_i)$ be a  5--dimensional
Seifert bundle, $L$ smooth.
Assume that  $H_1(L,\z)=0$  and $(S,\sum (1-\frac1{m_i})D_i)$
is a Calabi--Yau
orbifold. Then 
\begin{enumerate}
\item the minimal resolution of $S$ is a K3-surface, and 
\item $L$ is homeomorphic to the connected sum of at most
21 copies of $S^2\times S^3$.
\end{enumerate}
\end{cor}

Proof. By (\ref{seif.over.CY.prop}), $S$ is a simply
connected orbifold with 
trivial canonical class. In dimension 2
cyclic quotient
singularities with
trivial canonical class are Du Val.
Let $h:S'\to S$ be the minimal resolution.
Then $K_{S'}=h^*K_S$ and $H_1(S',\z)=0$.
Therefore $S'$ is a K3 surface and 
$\rank H^2(S,\q)\leq \rank H^2(S',\q)=22$, the maximum achieved
in the smooth case only.

(It may happen that $S'$ is a K3 surface but
$\pi_1^{orb}(S)\neq 1$. A famous example is the
Kummer surface, a quotient of an Abelian surface by $p\mapsto -p$.
In this case $\pi_1^{orb}(S)$ sits in an exact sequence
$$
0\to \z^4\to \pi_1^{orb}(S)\to \z/2\to 0.
$$
There are numerous similar cases, classified by
\cite{mukai}.)

Since there is no branch divisor, (\ref{simpconn.over.ratsurf.thm})
 shows that
$H_2(L,\z)$ is torsion free.
The  second
Stiefel--Whitney class is zero by (\ref{w2=0.lem}).
Therefore $L$ is homeomorphic to the connected sum of at most
21 copies of $S^2\times S^3$ by (\ref{smale.thm}).\qed
\medskip

The following lemma is essentially in    \cite{moroianu, bgn03c}.

\begin{lem}\label{w2=0.lem} Let $f:L\to (X,\Delta)$ be a Seifert bundle.
Assume that $H_1(L,\z)=0$ and that $K_X+\Delta$ is a
rational multiple of $c_1(L/X)$. Then $w_2(L)=0$.
\end{lem}

Proof. $L$ is an orientable hypersurface in the corresponding
Seifert $\c^*$-bundle $Y\to S$ and
$w_2(L)=c_1(Y)|_L\mod 2$. By \cite[41 and 16]{k-seif},
$$
c_1(Y)=-K_Y=f^*(K_X+\Delta)=c\cdot f^*c_1(Y/X)=0
$$ in rational cohomology.
Since  $H_1(L,\z)=0$, the second cohomology has no torsion
 thus $c_1(Y)=0$ and also $w_2(L)=0$.
\qed

\bibliography{refs}

\vskip1cm

\noindent Princeton University, Princeton NJ 08544-1000

\begin{verbatim}kollar@math.princeton.edu\end{verbatim}

\end{document}